

\input amstex.tex

\documentstyle{amsppt}
\magnification \magstephalf
\pagewidth{6.2in}
\pageheight{8.0in}
\NoRunningHeads
\def\makefootline{\baselineskip=49pt \line{\the\footline}}
\define\s{\sigma}
\define\la{\lambda}
\define\al{\alpha}

\define\ep{\varepsilon}

\define\th{\theta}

\define\pprime{{\prime\prime}}

\define\trace{\operatorname{trace}}

\define\intrinsic{{\operatorname{intrinsic}}}

\define\Hess{\operatorname{Hess}}

\define\End{\operatorname{End}}

\define\llangle{\langle\!\langle}
\define\rrangle{\rangle\!\rangle}
\define\Llangle{\left\langle\!\left\langle}
\define\Rrangle{\right\rangle\!\right\rangle}

\define\Supp{\text{Supp}}
  
\redefine\flat{{\text{flat}}}

\define\Hol{\operatorname{Hol}}

\def\today{\ifcase\month\or
    January\or February\or March\or April\or May\or June\or
    July\or August\or September\or October\or November\or December\fi
    \space\number\day, \number\year}

\topmatter
\title
HESSIANS OF SPECTRAL ZETA FUNCTIONS
\endtitle
\author  K. Okikiolu
\endauthor

\thanks  The author was supported  by the National
Science Foundation  \#DMS-9703329.  \endthanks \abstract Let $M$
be a compact manifold without boundary. Associated to a metric $g$
on $M$ there are various Laplace operators, for example the De
Rham Laplacian on $p$-forms and the conformal Laplacian  on
functions. For a general geometric differential operator of
Laplace type with eigenvalues $0\leq \la_1\leq\la_2\leq\dots$ we
consider the spectral zeta function $Z(s)=\sum_{\la_j\neq0}
\la_j^{-s}$. The modified zeta function $\Cal
Z(s)=\Gamma(s)Z(s)/\Gamma(s-n/2)$ is an entire function of $s$.
For a fixed value of $s$ we calculate the hessian of $\Cal Z(s)$
with  respect to the metric and show that it is given by a
pseudodifferential operator $T_s=U_s+V_s$ where $U_s$ is
polyhomogeneous of degree $n-2s$ and $V_s$ is polyhomogeneous of
degree $2$. The operators $U_s/\Gamma(n/2+1-s)$ and
$V_s/\Gamma(n/2+1-s)$ are entire in $s$. The symbol expansion of
$U_s$ is computable from the symbol of the Laplacian.  Our
analysis extends to describe the hessian of $(d/ds)^k \Cal Z(s)$
for any value of $k$.
\endabstract
\endtopmatter

\noindent{\bf 1. Introduction.}
\medskip

Let $M$ be an $n$-dimensional smooth compact manifold without
boundary. The tensor powers of $\Bbb R^n$  are $T^{p,q}(\Bbb
R^n)=\otimes^p\Bbb R^n\otimes\otimes^q{\Bbb R^n}^*$.
Suppose $U$ is a vector subspace  of $T^{p,q}(\Bbb R^n)$
which is invariant under the action of the general linear group
$GL(\Bbb R^n)$.  Writing $\Cal P$ for the principal frame bundle of $M$,
we obtain a general tensor bundle $E$ over $M$ by forming the
associated bundle
$$
 E =\Cal P\times_{GL(\Bbb R^n)} U.
 \tag 1.1
$$
For example we obtain in this way
$T^{p,q}(M)$, the bundle of  type-$(p,q)$ tensors
on $M$, and  $S^2M$, the bundle of symmetric $(0,2)$-tensors.
Let $\End(E)$ denote the bundle whose fibers are the endomorphisms of the
fibers of $E$, so $\End(E)$ can be identified with  $E\otimes E^*$.
Let $N$ denote the rank of $E$,
and let $C^\infty(E)$ denote the smooth sections of $E$.

For a
Riemannian metric $g$ on $M$, let $dV$ denote the canonical volume
element in the metric $g$, and $V$ the total volume of $M$. The
metric $g$ naturally gives a metric $\langle\,,\,\rangle_g$ on $E$
which gives an inner product $\llangle\,,\,\rrangle_g$ on
$C^\infty(E)$,
$$
\llangle\th,\omega\rrangle_g\ =\ \int_M \langle
\th,\omega\rangle_g\, dV,\qquad\qquad\qquad \th, \omega\in
C^\infty(E).
$$
The metric $g$ also gives rise to the Levi-Civita connection on
$M$ with curvature tensor $R$, and the connection on $E$
$$
\nabla=\nabla_g:C^\infty(T^{k} M\otimes E)\ \to\
C^\infty(T^{k+1}\otimes E),
$$
which is invariant under scaling, that is for $c>0$ and $\omega\in
T^k (M\otimes E)$ we have $\nabla_{c g}\omega=\nabla_g\omega$.

Now suppose that  for each Riemannian metric $g$ we have a differential
operator $F=F_g:C^\infty(E)\to C^\infty(E)$   which is
self-adjoint with respect to $\llangle\,,\,\rrangle_g$, non-negative,
scale invariant, and {\it  geometric of  Laplace type}, meaning that
$$
F\omega\ =\ V^{2/n}\left(\sum_{i=1}^n \nabla^2 (e_i,e_i,\omega) \,
+\,
\sum_{i,j,k,\ell=1}^nR(e_i,e_j,e_k,e_\ell)A(\th^i,\th^j,\th^k,\th^\ell,\omega)
\right)\ +\ B\omega \tag 1.2
$$
for $\omega\in C^\infty(E)$, where $B$ is a constant,
$e_1,\dots,e_n$, is a pointwise $g$-orthonormal base for $TM$
with dual coframe $\th^1,\dots,\th^n$, and
$A \in  C^\infty(T^{4,0}(M)\otimes \End(E))$  is obtained from some element
$$
\Cal A\ \in\ T^{4,0}(\Bbb R^n)\otimes U\otimes U^*
$$
which
is invariant under the natural action of $O(\Bbb R^n)$.
Indeed, if $\Cal O$ is the $g$-orthonormal
frame bundle of $M$, then
$$
T^{4,0}(M)\otimes E\otimes E^*\ \equiv\  \Cal O\ \times_{O(\Bbb R^n)}
\ (T^{4,0}(\Bbb R^n)\otimes U\otimes U^*)
$$
and $\Cal A$ naturally gives rise to  a section $A=A_g$ of
$T^{4,0}(M)\otimes E\otimes E^*$. Examples of $F$ include the de
Rham Laplacian on forms and the conformal Laplacian.

Write $\Re s$ for the real part of $s\in\Bbb C$.
Let $\la_1\leq \la_2\leq \dots $ be the eigenvalues of $F$ and set
$$
Z(s)\ =\ \sum_{j=1}^\infty \la_j^{-s}.
$$
Then $Z(s)$  converges for $\Re s>n/2$ and extends analytically to
a meromorphic function with possible simple poles at
$$
s\ =\ \cases n/2,n/2-1,\dots 1,\qquad & \text{for } n\text{ even}, \\
n/2,n/2-1,\dots \dots  & \text{for } n\text{ odd}     .\endcases
$$
In fact, if $F$ is invertible, then
 $$
 \Cal Z(s)\ :=\ \frac{\Gamma(s)Z(s)}{\Gamma(s-n/2)}
 $$
 is entire.  In general, when
$F$ has a kernel of stable dimension $d$, then
$$
 \Cal Z_{ g}(s)\ +\ \frac{d}{\Gamma(s-n/2)\ s}     \ \ \text{is
 entire.}
 \tag 1.3
 $$
and so for two metrics $g$ and $\tilde g$,
$\Cal Z_{\tilde g}(s)- \Cal Z_g(s)$ is
 entire.

 \remark{Remark} (1.3) does not hold for general
 polyhomogeneous pseudodifferential
 operators of stable dimensional kernel in odd dimensions,
 since $Z(s)$ may not vanish at the negative integers.
 See [Ok1] 1.2 for further details.
 \endremark

  With $s$  fixed,  differentiating $\Cal Z(s)$ twice with respect to
the metric one obtains the Hessian of $\Cal Z(s)$ which is a
symmetric bilinear form on the space $C^\infty(S^2M)$.  At the
metric $g$,
$$
\Hess \Cal Z(s)(h,h)\ =\ \frac{d^2}{d\s^2}\biggl|_{\s=0} \Cal
Z_{g+\s h}(s).
$$
Set $\Bbb N^+=\{1,2,\dots\}$.

 \proclaim{Theorem 1}  For $s\in \Bbb C$,  there exists a unique
symmetric pseudodifferential operator $T_s=T_s(F
):C^\infty(S^2M)\to C^\infty(S^2M)$ such that
$$
\Hess \Cal Z(s)(h,h)\ =\ \llangle h, T_s h\rrangle_g. \tag 1.4
$$
The operator $T_s$ is analytic in $s$.  For  $s\notin n/2+\Bbb
N^+$, there exist polyhomogeneous pseudodifferential operators
$U_s$ and $V_s$ of degrees $n-2s$ and $2$ respectively such that
$T_s=U_s+V_s$.  The operators $U_s$ and $V_s$ are meromorphic in
$s$ with simple poles located in $ n/2+\Bbb N^+$. (The poles of
$U_s$ and $V_s$ cancel  in the sum $T_s=U_s+V_s$, but the symbol
expansion of $T_s$ for $s\in n/2+\Bbb N^+$ involves logarithmic
terms.) For general $s$, the polyhomogeneous symbol expansion of
$U_s$ is computable from the complete symbol of the operator $F$.
In particular, there is a simple algorithm to compute the term
$u_s$ of homogeneity $n-2s$.
 Furthermore, we can differentiate in
$s$ to obtain
$$
\Hess (d/ds)^k \Cal Z(s)(h,h)\ =\ \llangle h, (d/ds)^k T_s
h\rrangle_g,\tag 1.5
$$
and the principal symbol of $(d/ds)^k U_s$ is equal to the leading
order  term of  $(d/ds)^k u_s$  (provided this does not vanish
identically).
\endproclaim

\proclaim{Theorem 2} The symbol $u_s$ from Theorem 1, can be
computed as follows.  Write $F'=(d/d\s)|_{\s=0} F(g+\s h)$ and
 let $x\in M$.  Take coordinates on $M$ which are orthonormal at
 the point $x$, and take a local trivialization of $E$ on a neighborhood of $x$ to
 obtain coordinates for $E$.
  Suppose that in these coordinates,  the
operator $F^\prime$ is given at the point $x$ by
$$
F^\prime\ =\ \sum_{\al,\beta,i,j} A_{\al,\beta}^{ij} \
(\partial^\al_w h_{ij}(x))\ \partial_w^\beta, \tag  1.6
$$
where  $\al$ and $\beta$ are multiindices and $A_{\al,\beta}^{ij}$
is an $N\times N$ real valued matrix.
 Set
$$
S=s-n/2,
$$
write $I$ for the identity operator on $E_x$, and set
$$
C(s)\ =\ \left(\frac1{4\pi}\right)^{n/2}
\frac{\Gamma(-S+1)^2}{\Gamma(-2S+2)}. \tag 1.7
$$
Then at $(x,\xi)\in T^*M$, the value of $u_s(x,\xi)\in\End
(S^2M)_x$ is given by
$$
(u_s(x,\xi)h)_{ij}\ =\  V^{(2s-n)/n}C(s)\sum\Sb |\al|+|\beta|=2\\
|\gamma|+|\delta|=2 \\ k,\ell\endSb
\,u_s(\partial^\al,\partial^\beta,\partial^\gamma,\partial^\delta,x,\xi)
\ \trace(A_{\al,\beta}^{ij}(x)\,A_{\gamma,\delta}^{k\ell}(x))\
h_{k\ell}, \tag 1.8
$$
where the terms
$u_s(\partial^\al,\partial^\beta,\partial^\gamma,\partial^\delta,x,\xi)$
are given as follows:
$$
u_s(\partial_j\partial_k, \, I,\, \partial_p\partial_q,\,I, \,
x,\xi) \ =\ 4(4S^2-1) \xi_j\xi_k\xi_p\xi_q\,|\xi|^{n-2s-4} .\tag a
$$
$$
u_s(\partial_j, \, \partial_k,\, \partial_p\partial_q,\,I,\,
x,\xi) \ =\   -2(4S^2-1) \xi_j\xi_k\xi_p\xi_q |\xi|^{n-2s-4}.\tag
b
$$
$$
u_s(\partial_j, \, \partial_k,\, \partial_p,\, \partial_q,\,
x,\xi) \ =\   (4S^2+2S-2)\xi_j\xi_k\xi_p\xi_q |\xi|^{n-2s-4}\ -\
 (2S-1)\delta_{kq}\xi_j\xi_p |\xi|^{n-2s-2}.\tag c
$$
$$
u_s(I, \, \partial_j\partial_k,\, \partial_p\partial_q,I,\, x,\xi)
\ =\  (4S^2-2S)\xi_j\xi_k\xi_p\xi_q |\xi|^{n-2s-4}\ +\
(2S-1)\delta_{jk}\xi_p\xi_q |\xi|^{n-2s-2}. \tag d
$$
$$
\multline u_s(I, \, \partial_j\partial_k,\, \partial_p,\,
\partial_q,\, x,\xi) \  =\ -(2S^2+S-1)\xi_j\xi_k\xi_p\xi_q
|\xi|^{n-2s-4}
 \\ +\   (S-1/2)(-\delta_{jk}\xi_p\xi_q+
 \delta_{jq}\xi_k\xi_p+\delta_{kq}\xi_j\xi_p) |\xi|^{n-2s-2}
\endmultline
\tag e
$$
$$
\align &u_s (I,  \,  \partial_j\partial_k, \,I,\,
\partial_p
\partial_q,\, x,\xi)
\  =  \tag f \\
 (S^2+ & S) \xi_j\xi_k\xi_p\xi_q|\xi|^{n-2s-4}\ +\
\frac12(S-1)\delta_{jk}\xi_p\xi_q|\xi|^{n-2s-2}
 \ +\ \frac12(S-1)\xi_j\xi_k\delta_{pq}|\xi|^{n-2s-2}\\
 -\   \frac{S}{2}& (\delta_{jp}\xi_k\xi_q+\delta_{kq}\xi_j\xi_p
 +\delta_{jq}\xi_k\xi_p+\delta_{kp}\xi_j\xi_q)|\xi|^{n-2s-2}\ +\
  \frac{1}{4}(\delta_{jk}\delta_{pq}+\delta_{jp}\delta_{kq}
  +\delta_{jq}\delta_{kp})|\xi|^{n-2s}.
\endalign
$$
\endproclaim

 We refer to the operator $T_s$ in Theorem 1, as {\it the operator corresponding to
$\Hess \Cal Z(s)$.} At a metric $g$ which is critical for $T_s$,
the eigenvalues of $T_s$ can determine whether  $g$ is a local
maximum, minimum or saddle point. Computation of  the principal
symbol of $T_s$ gives at least some information about the type of
critical point, for example if it changes sign, one can say for
sure that $g$ is a saddle point, while if it is positive, one can
generally determine that $g$ has finite index.
 The operators $U_s$, $V_s$ are
only well defined up to the addition (respectively subtraction) of
a smoothing operator which is analytic in $s$. They are however
canonically defined on the level of symbol expansions because in
general they have degrees which differ by a non-integer. It should
be noted that while the symbol expansion of $U_s$ is locally
computable, the expansion of $V_s$ generally involves global
quantities. Theorem 1 sheds light on the question of [GS]
concerning how many derivatives of an invariant must be computed
to obtain a local expression. We see that for low values of $s$,
it takes two derivatives of a zeta invariant to obtain an
expression whose ``leading order" is local. Indeed,  in general
$u_s$ is the principal symbol of $T_s$ when $\Re s<n/2-1$ and this
is local, while $v_s$ is the principal symbol when $\Re s>n/2-1$.
Of particular interest in this discussion is the zeta-regularized
determinant of $F$, which is defined to be $e^{-Z^\prime(0)}$. The
quantities  $u(0)$ and $(d/ds)|_{s=0} u_s$ usually control the
type of critical points possible for the determinant. The picture
may become complicated  at values  of $s$ for which $u_s$ or $v_s$
vanishes. At these values, the principal symbol of $T_s$ is given
by the highest degree non-vanishing  term in the symbol of $U_s$
or $V_s$, or if these terms have the same degree it is given by
the sum.  For example, in  [Ok2] we will
 actually compute the principal symbol of $V_s$ for the family of operators
 in (1.9).
 We will see cases when $V_s$ has
 order zero, and $u_s$ is the principal symbol of $T_s$ for $\Re
s<n/2$.    We also remark that one can formulate a slightly more
abstract but more basic version of Theorem 1, which does not
mention metrics or geometric operators. Indeed, fix a tensor
bundle over $M$, and consider the space of elliptic, non-negative
second order differential operators with isotropic principal
symbols.  Now consider the space of coefficients of such operators
and compute the Hessian of the zeta function with respect to
variation of the coefficients. Then
 this bilinear form on the space of coefficients is given by a
pseudodifferential operator, in a similar fashion to Theorem 1.
Theorem 1 then follows from the more abstract formulation by
expressing the coefficients of a Laplace type operator in terms of
the metric.

In this paper we prove Theorems 1 and 2.   Although  similar
results hold for geometric elliptic operators of  more general
form, we will not deal with greater generality here.  For us  the
main goal is to understand the behavior of zeta functions for
specific geometric operators. In two further papers [Ok2] and
[OW], we carry out specific calculations.  In [Ok2], we explicitly
compute $u_s$  for the family of Laplacians
$$
F\ =\ V^{2/n}\left(   \Delta \ +\ c_1 \mu S \right) \ +\ c_2 I ,
\tag 1.9
$$
where
$\Delta$ is the (positive) Laplace-Beltrami operator,
 $S$ is the scalar curvature,
$I$ is the identity operator,  $\mu=(n-2)/(4(n-1))$, and $c_1,c_2$
are arbitrary constants which render $F$ non-negative. This family
contains both the Laplace-Beltrami operator and the conformal
Laplacian ($c_1=1$, $c_2=0$), and we are in fact able  to compute
the principal symbol of $V_s$ which is non-local. We go on to make
conclusions about the critical points of $Z(s)$ for all values of
$s$. For large values of $s$, $Z(s)$ should be dominated by the
lowest eigenvalues, and one sees consequenses of this in the
results. In joint work with C. Wang, [OW], we compute $u_s$ for
the de Rham Laplacian.  Theorem 2 (a)-(f) gives the  symbol $u_s$
needed for the applications [Ok2] and [OW]. The general formula
for $u_s$ is given in Lemma 5.1.  We do not give here explicit
formulas for the complete symbol expansions of $U_s$ and $V_s$,
although the symbol expansion of $U_s$ can be obtained from (5.10)
or Lemmas 6.1, 6.5, and the symbol expansion of $V_s$ has
contributions from three  terms, firstly (4.5), secondly (4.16)
whose symbol expansion can be computed from  (4.18), and thirdly
(4.8) whose symbol expansion is computed in   Lemmas 6.1, 6.5.

This work is
motivated by a growing list of results [OPS1], [BCY], [CY], [Br],
 [Ri], [Ok1], [Mo1], [Mo2], where the authors take $F$ to be
the conformal Laplacian or Laplace-Beltrami operator, and identify
local or global maximal or minimal metrics for the functional
$\log\det F$ or $Z(s)$, either within a conformal class, or among
all smooth metrics. Papers [OPS2], [OPS3] and [CGY] contain
applications of some of those results. It emerges that in some
situations only maxima occur  while in others only minima occur.
In [Ok2] and [OW] we will address   the questions of when it is
possible for maxima to occur, when it is possible for minima to
occur, when neither can happen, and when one might hope to find a
convexity argument to prove global results.
\medskip
\definition{Acknowledgements}
The author would like to thank the referees for their careful and
extensive comments on earlier versions of this manuscript.  In
addition, the author would also like to thank the following
institutions for support: ESI, IAS, IHES, MSRI, Princeton
University.
\enddefinition

We will now indicate briefly what is involved in the proof of
Theorems 1 and 2. One can use the parametrix for the heat kernel
or the resolvant of the Laplacian to prove the result. Technically
the work involved is  comparable for the two methods, and each
method has its advantages. The heat kernel method has the
advantage of yeilding completely explicit formulas for the symbol
expansion of $U_s$, and so we take this approach. We fix a metric
$g$, and take a smooth one-parameter deformation $\tilde g(\s)$
with $\tilde g(0)=g$. Using the Melin transform to write $\Cal
Z_{\tilde g(\s)}(s)$ in terms of the heat kernel, we can
differentiate to compute the second derivative of $\Cal Z_{\tilde
g}(s)$ with respect to $\s$ at $\s=0$, see Lemma 3.1. The main
term we need to study is the one giving rise to the operator
$U_s$, which is the last term on the right hand side of (3.1).
Writing $F^\prime= (\partial F_{\tilde g(\s)}/{\partial\s})(0)$,
it equals
$$
 \frac{1 }{\Gamma(s-n/2)}
\trace  \iint_{u+v<1}\  (u+v)^{s}
 F^\prime \ e^{-u F }\    F^\prime   \  e^{-v F}\, du
 dv. \tag 1.10
$$
Before proceeding, we give a list of the steps in the proof to
help the reader navigate.  The steps in this guide will usually
be referenced when they occur in the proofs that follow.

\proclaim{Table 1.1: Steps of the Proof} This table lists the
steps which are needed to go from the expression (1.10)
 to the symbol expansions of the operators
$U_s$ and $V_s$.  They are not listed here in the order they occur
in the proof,   and some of the steps need to be repeated  more
than once in the proof. Steps 1-8 are entirely trivial while
Steps 9--16 represent the analytic details, with Steps 13--16
being the more significant or subtle. On the other hand, one has
to pay a little attention to steps 7, 8 and 9, because they change
the nature of the expression being considered.
 \roster

\item""{\bf Step 1}:  The steps just carried out produced a sum of terms.
Consider each term  separately.

\item""{\bf Step 2}:  At a previous stage, the expression was
decomposed as a sum and each term was considered separately.
Reverse this by summing up the contributions from each term.

\item"" {\bf Step 3}: Introduce a partition of unity to
reduce to coordinate charts over which $E$ is trivial.

\item""{\bf Step 4}: Work in geodesic normal
coordinates.

 \item""{\bf Step 5}:  Write $h=\tilde g'(0)$ in components.

\item""{\bf Step 6}: Write $F'$ in components.

\item""{\bf Step 7}: The expression has the form ``$\iint
((u+v)^s\ \ \overline{\hphantom{aaa}}\ \ dudv$". Restrict
attention to the term ``$\overline{\hphantom{aaa}}$" in the
integrand.

\item""{\bf Step 8}: Reverse Step 8 by reinstating the factor
$(u+v)^s$ and the $uv$ integration.

\item""{\bf Step 9}: Identify an operator $\tilde T$ (by its
kernel) so that the expression has the form $\llangle h,\tilde T
h\rrangle$.

\item""{\bf Step 10}: Write the trace of an operator as the
integral of the kernel over the diagonal.

\item""{\bf Step 11}: Make the standard small time expansion of
the heat kernel.

\item""{\bf Step 12}: Make a Taylor expansion of each term in a
given variable.

\item""{\bf Step 13}: Collect together terms of the same ``order".

\item""{\bf Step 14}: Exchange $uv$ integration in the expression
with another limiting process such as taking the operator trace, a
Fourier transform or integration in other variables.

\item""{\bf Step 15}: Compute the symbol of the operator from its
kernel by taking the Fourier transform.

\item""{\bf Step 16}: The symbol has the form $\iint ((u+v)^s\ \
\overline{\hphantom{aaa}}\ \ dudv$". Carry out the $uv$
integration to see the asymptotic behavior in the cotangent
variable $\xi$.

\endroster
\endproclaim

\noindent{\bf The simplified scalar case}:  We wish to write
(1.10) in terms of $h=(d\tilde g/d\s)(0)$.
 Let us first pretend for the sake of clarity that
$F^\prime$ is just multiplication by the smooth function
$\phi=h_{ij}A^{ij}$, where $A^{ij}$ is a smooth tensor field. In
doing so we will be ignoring some technicalities which arise when
$A$ is a differential operator, but let us concentrate here on the
main point. This eliminates the need for {\bf Steps 5, 6}. We now
wish to compute
$$
 \frac{1 }{\Gamma(s-n/2)}
\trace  \iint_{u+v<1}\  (u+v)^{s}
 \phi \ e^{-u  F }\    \phi   \  e^{-v F}\, du
 dv.\tag 1.11
$$
{\bf Steps 14, 10}: Writing $K(e^{-tF},x,y)$ for the integral
kernel of $e^{-tF}$, this is equal to
$$
\frac{1 }{\Gamma(s-n/2)}   \iint_{u+v<1}\  (u+v)^{s}
 \iint_{M\times M} \phi(x) \ K(e^{-u  F },x,y)\    \phi(y)
 \ K( e^{-v F},y,x)\, dV(x) dV(y)\, du
 dv.\tag 1.12
$$
{\bf Steps 14, 9}: We see that this equals
$$
\llangle \phi, \tilde T_s\phi\rrangle\tag 1.13
$$
where $\tilde T_s$ is the operator on functions with
$$
(\tilde T_s\phi)(x)\ =\ \int_M K_s(x,y) \phi(y)\, dV(y),\tag 1.14
$$
where
$$
K_s(x,y)\ =\ \frac{1 }{\Gamma(s-n/2)}   \iint_{u+v<1}\  (u+v)^{s}
  \ K(e^{-u  F },x,y)  \ K( e^{-v F},y,x)\, du
 dv.\tag 1.15
$$
To express this operator as a psuedodifferential operator, we
compute its symbol by taking the Fourier transform.  Here, we
focus just on the leading order term in this symbol.

\noindent{\bf Steps 3, 1, 4}:  Fix a point $x\in M$ and work in
normal coordinates $w$ of $y$ about $x$.

\noindent{\bf Steps 7, 11, 13}: For this discussion, let us assume
$V=1$, and restrict attention to the leading order term of the
heat kernel $K(e^{-tF}, x,y)$, which   a constant multiple of, $
t^{-n/2} e^{-|w|^2/4t}$. Approximately then we have
$$
K(e^{-u  F },x,y)  \ K( e^{-v F},y,x)\ \approx\ C( uv)^{-n/2}
e^{-|w|^2(1/4u+1/4v)}.\tag 1.16
$$
{\bf Step 15}: Since $\sqrt{g(x)}=1$, up to leading order, the
symbol of the right hand side of (1.16) is obtained by taking the
Fourier transform in $w$. This yields a constant multiple of,
$$
(u+v)^{-n/2} e^{-|\xi|^2 uv/(u+v)}.\tag 1.17
$$
{\bf Steps 14, 8}: From this, we see that at $x$, the symbol of
$\tilde T_s$, whose kernel is given by (1.15), is approximately
$$
\frac{1 }{\Gamma(s-n/2)}\iint_{u+v<1}\  (u+v)^{s-n/2}
 e^{-|\xi|^2 uv/(u+v)}\, du dv.\tag 1.18
 $$
{\bf Step 16}: Substituting $\tau=u/(u+v)$ and $T=|\xi|^2(u+v)$,
this is
$$
 \frac{|\xi|^{n-2s-4} }{\Gamma(s-n/2)}\int_{T<|\xi|^2}\
\int_{\tau=0}^1 T^{s-n/2+1}
 e^{-T \tau(1-\tau)}\, d\tau dT. \tag 1.19
 $$
 If  we  replace the interval of
 integration by $0<T<\infty$, we
 get
 $$
 \frac{|\xi|^{n-2s-4 } }{\Gamma(s-n/2)}\int_{T=0}^\infty\
\int_{\tau=0}^1 T^{s-n/2+1}
 e^{-T \tau(1-\tau)}\, d\tau dT\ =\
 \frac{\Gamma(n/2-s-1)\,\Gamma(n/2-s+1)\
 |\xi|^{n-2s-4}}{\Gamma(n-2s-2)}. \tag 1.20
 $$
 This gives the principal symbol of $U_s$. The reason it
 has degree $n-2s-4$ instead of $n-2s$ as stated in Theorem 1, is because
 we have assumed that $F^\prime$ is just multiplication by the
 function
 $\phi=h_{ij}A^{ij}$, whereas in general it may involve derivatives up to order
 $2$.  Since $F^\prime$ occurs twice in (1.10), this accounts for the
  four missing derivatives.  To compute (1.19) we need
to subtract from (1.20) the part of the integral on the left hand
side  over the interval $ |\xi|^2<T<\infty$. But by steepest
descent, there is an expansion as $T\to\infty$ of the form
 $$
\int_{\tau=0}^1
 e^{-T \tau(1-\tau)}\, d\tau\ \sim\ 2T^{-1}\ +\ C_2 T^{-2}\ +\ C_3
T^{-3}\ +\ \dots.\tag 1.21
 $$
Hence as $|\xi|\to\infty$,
 $$
 \multline
\frac{|\xi|^{n-2s-4 } }{\Gamma(s-n/2)}\int_{T>|\xi|^2}T^{s-n/2+1}\
\int_{\tau=0}^1
 e^{-T \tau(1-\tau)}\, d\tau dT\\
  \sim\  \frac{-1}{\Gamma(s-n/2)}\left(\frac{2|\xi|^{-2}}{s-n/2+1}\
  +\
\frac{C_2 |\xi|^{-4}}{s-n/2}\ +\ \frac{C_3 |\xi|^{-6}}{s-n/2-1}\
+\ \dots\right).
\endmultline
\tag 1.22
 $$
This is a symbol of order $-2$.  The  poles, which occur at
$s=n/2+1,\,n/2+2,\dots$, cancel with the poles of (1.20) and
contribute logarithmic terms to the symbol of $\tilde T_s$ at
these values of $s$.
\medskip

\noindent{\bf The general case}: Now we modify this basic
calculation to obtain Theorem 2 on the tensor bundle $E$.

\noindent{\bf Step 5}:  Suppose that $h$ is supported on a small
neighborhood $\Omega$ of a point $x$ on which we take coordinates.

\noindent{\bf Step 6}: Write $F^\prime$ in components on $\Omega$
as
$$
F'\ =\ \sum_{\al,\beta,i,j} A_{\al,\beta}^{ij} \ (\partial^\al_w
h_{ij})\ \partial_w^\beta,
$$
 so that (1.10) becomes
$$
\sum\Sb\al,\beta,\gamma,\delta\\i,j,k,\ell\endSb  \frac{1
}{\Gamma(s-n/2)}\trace \iint_{u+v<1}\ (u+v)^{s} A_{\al\beta}^{ij}\
(\partial^\al h_{ij})\,\partial^\beta \ e^{-u F }\
A_{\gamma\delta}^{k\ell}\ (\partial^\gamma h_{k\ell})\,
\partial^\delta   \ e^{-v F}\, du
 dv. \tag 1.23
$$
In fact, because we are only interested in the leading order term
in the symbol, we can restrict the sum in (1.23) to those
multiindices with $|\al+\beta|=|\gamma+\delta|=2$.

\noindent{\bf Steps 1, 14, 10}:  We will analyze each term in the
sum (1.23) separately. We write
$$
\multline \frac{1 }{\Gamma(s-n/2)}\trace \iint_{u+v<1}\ (u+v)^{s}
A_{\al\beta}^{ij}\ (\partial^\al h_{ij})\,\partial^\beta \ e^{-u F
}\ A_{\gamma\delta}^{k\ell}\ (\partial^\gamma h_{k\ell})\,
\partial^\delta   \ e^{-v F }\, du
 dv\\
 =\ \frac{(-1)^{|\al+\gamma|} }{\Gamma(s-n/2)}   \iint_{u+v<1}\  (u+v)^{s}
 \iint_{M\times M}(\partial^\al h_{ij})(x)\ E(u,v,x,y)
 \ (\partial^\gamma h_{k\ell})(y)
 \, dV(x)dV(y)\, du
 dv,
 \endmultline
 \tag 1.24
 $$
 where
 $$
 E(u,v,x,y)\ =\ \trace \ A_{\al\beta}^{ij}(x)\ (\partial_x^\beta K(e^{-u  F  },x,y))\
 A_{\gamma\delta}^{k\ell}(y)\ (\partial_y^\delta K( e^{-v F },y,x)).
 \tag 1.25
 $$
 {\bf Steps 14, 9}: This has the form
 $$
 \llangle h_{ij}, T_{s,ij}^{k\ell} h_{k\ell}\rrangle,
 \tag 1.26
 $$
 where $T_{s,ij}^{k\ell}$ is the operator on scalar valued
 functions whose kernel to leading order is
$$
 \frac{(-1)^{|\al+\gamma|} }{\Gamma(s-n/2)}   \iint_{u+v<1}\  (u+v)^{s}
  \partial_x^\al\partial_y^\gamma E(u,v,x,y)
 \ \, dudv. \tag 1.27
$$
(We integrated by parts to obtain this.  Derivatives falling on
the volume factors $\sqrt{g(x)}$ and $\sqrt{g(y)}$ yield lower
order terms.)

 \noindent
{\bf Steps 11, 13}: To compute the principal symbol of the kernel
in (1.27), we restrict attention the leading order term of the
heat kernel $K(e^{-tF }, x,y)$, which is a constant multiple of, $
t^{-n/2} e^{-|d(x,y)|^2/4t}\,I$, where $d(x,y)$ is the distance
from $x$ to $y$.  This enables us to replace (1.27) by
$$
\multline
 \frac{(-1)^{|\beta+\gamma|}\trace A_{\al\beta}^{ij}(x)\
 A_{\gamma\delta}^{k\ell}(x)  }{\Gamma(s-n/2)}  \\
 \times\  \iint_{u+v<1}\,(u+v)^{s}\  (uv)^{-n/2} \,
  \left(\partial_y^\al\partial_y^\gamma \left(
 (\partial_y^\beta e^{-|d(x,y)|^2/4u})(\partial_y^\delta e^{-|d(x,y)|^2/4v})\right)
 \right) \ \, du dv.
 \endmultline\tag 1.28
$$
Here, we are using the fact that to leading order,
$$
\partial_x^\al e^{-d^2(x,y)/4u}\bigl|_{x=0}\ \sim\ (-1)^{|\al|}
\ \partial_y^\al e^{-d^2(x,y)/4u}.
$$

\noindent{\bf Step 4}: At this point we fix the point $x$ and work
in geodesic normal  coordinates on $\Omega$ centered at $x$. (The
coordinates used in (1.23) to reduce to a scalar problem may be
arbitrary.)

\noindent {\bf Steps  14, 15, 16}: Taking the  Fourier transform
of (1.28) in the normal coordinates of $y$ as in (1.17)--(1.20),
yields the degree $n-2s$ term in the symbol of $ T_{s,ij}^{k\ell}$
of (1.26). Carrying out the calculation for explicit choices of
$(\al,\beta,\gamma,\delta)$ gives rise to the quantities  in
Theorem 2, (a)-(f).
\medskip

The rest of the paper is organized as follows. In Section 3, we
derive variation formulas for $\Cal Z(s)$ to write  the Hessian of
$\Cal Z(s)$ in terms of derivatives of $F$. In Section 4, we  use
basic estimates on the heat kernel to reduce the problem of
studying the Hessian to that of studying (1.10) (equivalently
(4.19)). In order to analyze the  operator associated to (1.10),
we will need to make a double decomposition of the integral
kernel, one stage to extract the polyhomogeneous expansion of the
symbol of $U_s$, and the other to extract the polyhomogeneous
expansion of $V_s$. In Section 5, we extend the above sketch to
obtain the expansion of the symbol of $U_s$ and use this to deduce
Theorem 2. However, we leave the rigorous estimates on errors to
Section 6 where we prove Theorem 1 by carrying out the double
decomposition in detail. Formula (1.5) then requires no proof
because it follows directly from (1.4), the fact that $T_s$ is an
analytic family of operators, and the standard fact that if the
kernel of $F$ has stable dimension then $\Cal Z_g(s)$ is smooth on
$\Cal M\times \Bbb C$.
\medskip\medskip


\noindent {\bf 2. Background and conventions.} Before commencing
with the proof of Theorem 1, we discuss in more detail the Hessian
of the function $\Cal Z(s)$ and we fix the conventions which will
be assumed concerning polyhomogeneous pseudodifferential operators
($\Psi$DOs). Let $\Cal M $ denote the space of Riemannian metrics
on $M$. Define the Sobolev  space $H^r(S^2M)$ to be the completion
of $C^\infty(S^2M)$ in the norm $\|\ \|_r$ defined by
$$
\|h\|_r\ =\ \sum_{s\leq r}\llangle \nabla^s h,\nabla^s
h\rrangle_g^{1/2},
$$
where $g$ is a fixed but arbitrary smooth metric. The Sobolev
imbedding theorem shows that for $r>n/2$, $H^r(S^2M)$  is
contained in the space $C(S^2M)$ of continuous sections of $S^2M$.
The smooth topology on $\Cal M$ is the Frechet space topology
coming from the collection of all the norms $\|\ \|_r$.  For
$r>n/2$, let $\Cal M^r$ denote the closure of the space of smooth
metrics on $M$ in the space $H^r(S^2M)$. For a compact subset
$U\subset \Bbb C$, let $\Hol(U)$ be the space of holomorphic
functions on $U$ with the supremum norm.

 \proclaim{Proposition} Suppose $m,r\in \Bbb N$ with $r>m+n/2+1$,
 and
 $U\subset\{s\in\Bbb C:\Re s> m+n/2+1-r\}$.  Then the map
 $g\to\Cal Z_{F _g}$
 is an $m$-times continuously differentiable  function from $\Cal M^r$ to $\Hol(U)$.
 \endproclaim

  The proof is given in [Ok3].  Using
 it, we obtain the Hessian of $(d/ds)^k\Cal Z(s)$ as a bilinear form on
 $H^r(S^2M)$, provided $r>n/2+3-\min\{(0,\Re s)\}$.  This restricts to a bilinear form
on $C^\infty(S^2M)$. We should
 remark that the Hessian defined above Theorem 1 is computed using
 the local linear space structure of $\Cal M$, that is using the
 flat connection on $\Cal M$.  It is possible to take other connections
 on $\Cal M$.
For a sufficiently regular function $F:\Cal M\to \Bbb C$ and a
smooth curve $\tilde g(\s)$ in $\Cal M$, the derivative $DF$ is
the linear form on $C^\infty(S^2M)$ defined by
$$
DF\left(\frac{\partial {\tilde g}}{\partial\s}\right)\ =\
\frac{\partial  F({\tilde g}(\s)) }{\partial \s}.
$$
The Hessian of  $F$ relative to a connection $\nabla$ on $\Cal M$
is given at $g\in\Cal M$ by
$$
\Hess F\left(\frac{\partial {\tilde g}}{\partial\s},
\frac{\partial {\tilde g}}{\partial\s} \right) \ =\
\frac{\partial^2 F(\tilde g(\s))}{\partial \s^2} \ -\ D
F\left(\nabla\left(\frac{\partial \tilde
g}{\partial\s},\frac{\partial \tilde g}{\partial
\s}\right)\right),\tag 2.1
$$

 We will prove Theorem 1 when the Hessian  is defined relative  to
 the flat or intrinsic connection on $\Cal M$.
The intrinsic connection on $\Cal M $ is the Levi-Civita
connection coming from the intrinsic metric defined as follows.
For $u,v\in C^\infty(S^2M)$, at the point $g\in \Cal M$,
$$
G_{\text{intrinsic}} ( u,v)\ =\  \llangle u, v \rrangle_g.
$$
Then
$$
 \left(\nabla\left(\frac{\partial \tilde g}{\partial
\s},\frac{\partial \tilde g}{\partial \s}\right)\right)_{ij}\ =\
 \frac{ \partial^2\tilde g_{ij}}{\partial\s^2} \
-\  \tilde g^{k\ell}\frac{\partial\tilde g_{ik}}{\partial\s}
\frac{\partial\tilde g_{\ell j}}{\partial\s} +\  \frac12\tilde
g^{k\ell}\frac{\partial\tilde g_{k\ell}}{\partial\s}
\frac{\partial\tilde g_{i j}}{\partial\s}
 -\ \frac14\tilde g_{ij}\tilde g^{k\ell}\tilde g^{pq}
 \frac{\partial\tilde g_{\ell p}}{\partial\s}
 \frac{\partial\tilde g_{kq}}{\partial\s}. \tag 2.2
$$
\medskip

Now we give a summary of some facts we will need about
polyhomogeneous $\Psi$DOs. For further details the
reader is referred to [Tr] and [Se]. If $U\subset \Bbb R^n$ is open
and $m\in \Bbb R$,
a function $q\in C^\infty(U\times\Bbb R^n,T^{1,1}(\Bbb R^N))$ is an
{\it $N\times N$-matrix valued
symbol of order $m$ (and type $(0,1)$)}, if for all multiindices
$\al,\beta$ and compact sets
$K\subset U$,
$$
|\partial_x^\al\partial_\xi^\beta q(x,\xi)|\ \leq\ C(1+|\xi|)^{m-|\beta|},
$$
where $C=C(\al,\beta,K)$ and $|\,\cdot\,|$ is a norm on
$T^{1,1}(\Bbb R^N)$. The symbol $q$ has order $-\infty$ if it has
order $m$ for all $m\in \Bbb R$. The symbol $q$ is  {\it
polyhomogeneous symbol of degree $\mu\in \Bbb C$} if there exists
a sequence of functions $q_j\in C^\infty(\Omega\times \{\xi\in\Bbb
R^n:\xi\neq0\},T^{1,1}(\Bbb R^N))$, $j=0,1,2,\dots$, which are
homogeneous;
$$
q_j(x,r\xi)\ =\ r^{\mu-j} q_j(x,\xi),\qquad\qquad\qquad r>0,
$$
such that
$$
q(x,\xi)\ \sim\ q_0(x,\xi)\ +\ q_1(x,\xi)\ +\ \dots,
\qquad\qquad\qquad\text{as } |\xi|\to\infty.
$$
By this we mean that if we introduce a cut-off function
$\phi\in C^\infty([0,\infty))$ supported on $[1/2,\infty)$
with $\phi=1$ on $[1,\infty)$, then
$$
q(x,\xi)\ -\  \phi(|\xi|)\sum_{j=0}^{J-1} q_j(x,\xi)
\ \text{ is a symbol of order }\Re \mu-J.
$$
(Where $\Re \mu$ is the real part of $\mu$.)

Now suppose $\pi:E\to\Omega$ is a rank $N$ vector bundle and
$Q:C^\infty(E)\to C^\infty(E)$ is a linear operator. Then $Q$ is a
$\Psi$DO of order $m$ if whenever we can trivialize $E$ over a
coordinate neighborhood thus identifying it with
 $(U\subset \Bbb R^n)\times \Bbb R^N$, then
there is a symbol $q$ of order $m$ such that
$$
Qh(x)\ =\ \frac1{(2\pi)^n}\int_{\Bbb R^n} \int_M
e^{i(x-y)\cdot\xi} q(x,\xi) h(y)\, dy d\xi. \tag 2.3
$$
for $x$ in $U$ and $h$ supported on a compact subset of $U$. The
principal symbol is then a section of the pull back of the bundle
$E\otimes E^*$ to the cotangent bundle of $M$.
 $Q$ is polyhomogeneous of degree $\mu$ if and only if for all
 trivializations of $E$, the symbol $q$ of $Q$ is polyhomogeneous
 of degree $\mu$.
We remark that from this it follows that there is a smooth kernel
$k$ on $M\times M\setminus \{(x,x)\}$ such that
$$
Qh(x)\ =\ \int_M k(x,y)h(y)\,dy
$$
whenever $x$ is not in the support of $h$.  Conversely, if the
operator $Q:C^\infty(E)\to C^\infty(E)$ is given by a smooth
kernel away from the diagonal then it is a $\Psi$DO of order $m$
provided (2.3) holds on  coordinate charts of small diameter. If
the operator $Q$ is given by a kernel which is smooth everywhere,
it is called a {\it smoothing operator}.

When the manifold $M$ is equipped with  a Riemannian metric $g$,
it is convenient to use the density $dV$ to define the kernel
$K(Q,x,y)$ of the operator $Q$, so
$$
Qh(x)\ =\ \int_M K(Q,x,y)h(y)\,dV(y).
$$
With this convention, the symbol $q(x,\xi)$ of $Q$ in a given
coordinate system is obtained by taking  the Fourier transform of
$K(Q,x,y)dV(y)$ :
$$
q(x,\xi)\ =\  \int K(Q,x,y) \chi(y)\, e^{i(y-x)\cdot \xi}\, dV(y).
$$
Here, $\chi$ is a cut off function equal to one on a neighborhood
of $x$ and supported in the coordinate chart. The symbol $q$ is
well defined up to a symbol of order $-\infty$ arising from the
choice of $\chi$.  In particular, if $Q$ is polyhomogeneous symbol
expansion is well defined at $x$ in the given coordinate system.
It will generally be convenient to compute the {\it
metric-modified symbol}
$$
q(x,\xi)\ =\  \int K(Q,x,y) \chi(y)\, e^{i(y-x)\cdot \xi}\, dy
\tag 2.4
$$
rather than the true symbol.  The two are of course very simply
related.

Suppose now that  we choose a smooth function
$w:\Omega\times\Omega\to\Bbb R^n$ such that for $x$ fixed,
$w(x,y)$ defines geodesic normal coordinates of $y$ centered at
$x$. Then  $(x,y)\to (x,w)$ defines a smooth diffeomorphism  from
a neighborhood of the diagonal in $\Omega\times\Omega$ into
$\Omega\times \Bbb R^n$.  It follows from [H\"{o}] that (2.3) can
be replaced by
$$
Qh(x)\ =\ \frac1{(2\pi)^n}\int_{\Bbb R^n} \int_M
 e^{iw\cdot\xi} \tilde q(x,\xi) h(y(x,w))\, dw d\xi,
$$
whenever there is a set of small diameter in $M$ which contains
both the support of  $h$  and the point $x$. Here, $\tilde q$ is a
symbol of order $m$ if we express $x$ in coordinates.

If $Q_s:C^\infty(E)\to C^\infty(E)$ is a family of $\Psi$DOs
depending on a parameter $s\in \Cal O\subset \Bbb C$,
then we say that $Q_s$ is a {\it analytic
family of  polyhomogeneous  $\Psi$DOs } if $Q_s$ has degree
$\mu(s)$ where $\mu$ is
analytic on $\Cal O$, and on any local trivialization of  $E$,
$Q_s$ has symbol
$q(s,x,\xi)$ which  is smooth in $(s,x,\xi)$  analytic in $s$,
and satisfies the estimates
$$
|\partial_x^\al\partial_\xi^\beta q(s,x,\xi)|
\ \leq\ C(1+|\xi|)^{\Re \mu(s)-|\beta|},
$$
where $\Re \mu(s)$ is the real part of $\mu(s)$ and $C$,
which depends on $\al,\beta$,  is uniform for
$x$ and $s$ in compact sets.  This implies that the homogeneous
terms $q_j(s,x,\xi)$ in the symbol expansion
of $Q_s$ are analytic in $s$, that the kernel defining $Q_s$ away
from the diagonal is analytic in $s$,
and that if $A:C^\infty(E)\to C^\infty(E)$ is a positive elliptic
second order differential operator,
then $Q_s A^{-\mu(s)/2}$ is a analytic family of bounded operators on $L^2(E)$.

\medskip

\medskip
\medskip
\medskip


\noindent{\bf 3. Variation formulas.}
\medskip
\medskip

 Let $g$ be a metric
on $M$, let $h\in C^\infty(S^2M)$, and  set
$$
\tilde g(\s )\ =\ g\ +\ \s \, h.
$$
We use primes to denote differentiation with respect to $\s$, so for example
$$
F^\prime=\frac{d}{d\s }F_{\tilde g(\s
)},\qquad\qquad F^\pprime=\frac{d^2}{d\s^2}F_{\tilde g(\s )}.
$$
In this section we will prove the following Lemma expressing  the
second variation of $Z(s)$ in terms of $F$, $F^\prime$ and
$F^\pprime$. \proclaim{Lemma 3.1}
 $$
\align \frac{d^2 \Cal Z(s)}{d \s ^2} \ =\ &\trace F^\pprime L(s) \
-\  \frac{2}{\Gamma(s-n/2)}\trace\left( F^\prime \Pi F^\prime
F^{-1}
 \int_{1}^\infty\  t^{s}\, (e^{-t F}-\Pi)\, dt \right)\\
 &\ +\  \frac{1 }{\Gamma(s-n/2)} \trace  \iint_{u+v>1}\  (u+v)^{s}
 F^\prime  \ (e^{-u  F }-\Pi)\    F^\prime   \ ( e^{-v F}-\Pi)\, du dv
\\ &\qquad\qquad\qquad\qquad\qquad\qquad +\  \frac{1 }{\Gamma(s-n/2)}
\trace  \iint_{u+v<1}\  (u+v)^{s}
 F^\prime  \ e^{-u  F }\    F^\prime   \  e^{-v F}\, du dv .
 \tag 3.1
\endalign
$$
\endproclaim

We start by expressing the first variation  in terms of $F$ and $F^\prime$.
In the case when the kernel of $F$ is trivial,
$$
\Cal Z(s)\ =\ \frac{\Gamma(s)}{\Gamma(s-n/2)}\trace F^{-s}, \tag
3.2
$$
and
$$
D  \Cal Z(s) (h)\ =\ \frac{d }{d\s}\Cal Z(s)\ =\
\frac{-\Gamma(s+1)}{\Gamma(s-n/2)} \trace F^\prime F^{-s-1}.
$$
When $F$ has a non-trivial kernel whose dimension does not vary with $\s$,
let $\Pi$ denote the projection onto the kernel, and write
 $$
F^{-s}\ =\  (F+\Pi)^{-s}\ -\ \Pi,
$$
so (3.2) still holds. Using $\Pi^2=\Pi$ we get
 $$
 \Pi^\prime \,\Pi\ +\ \Pi\, \Pi^\prime\ =\ \Pi^\prime,
 $$
 and using $F\Pi=0$ we get
 $$
 F^\prime \Pi\ =\ -F\,\Pi^\prime
 $$
 so
 $$
 \Pi F^\prime \Pi\ =\ 0
 $$
 and
 $$
 \trace F^\prime \Pi\ =\ 0.
 $$
 Hence
$$
\trace \Pi^\prime(F+\Pi)^{-s-1}\ =\ 2\trace \Pi^\prime\Pi(F+\Pi)^{-s-1}
\ =\ 2\trace \Pi^\prime\Pi\ =\ \trace\Pi^\prime\ =\ 0,
$$
and
$$
\multline D \Cal Z(s) (h)\ =\ \frac{-\Gamma(s+1)}{\Gamma(s-n/2)}
\trace (F+\Pi)^\prime (F+\Pi)^{-s-1} \ =\
\frac{-\Gamma(s+1)}{\Gamma(s-n/2)} \trace F^\prime (F+\Pi)^{-s-1}\\
=\ \frac{-\Gamma(s+1)}{\Gamma(s-n/2)} \trace F^\prime F^{-s-1}.
\endmultline
$$
It will be convenient to write this in the form
$$
DZ(s)(h)\ =\ \trace F^\prime L(s) \tag 3.3
$$
where
$$
L(s)\ =\ \frac{-1}{\Gamma(s-n/2)} \left( \Gamma(s+1)F^{-s-1}\ +\
\frac1{s+1}\Pi\right). \tag 3.4
$$

Now we compute the second variation of $\Cal Z(s)$ in terms of $F$, $F^\prime$,
and $F^\pprime$.
In the case when the kernel of $F$ is trivial,
$$
\frac{d^2 \Cal Z(s)}{d\s^2}\ =\ \frac{-\Gamma(s+1)}{\Gamma(s-n/2)}\trace F^\pprime
F^{-s-1}
 \ -\  \frac{\Gamma(s+1)}{\Gamma(s-n/2)}\trace F^\prime \frac{ d\tilde F^{-s-1}}{d\s }.
$$
Using the Melin transform, we write
$$
\tilde F^{-s-1}\ =\ \frac1{\Gamma(s+1)}   \int_0^\infty t^{s}
e^{-t\tilde F}\, dt.
$$
By Duhamel's principle,
$$
\frac{\partial e^{-t\tilde F} }{\partial \s } \ =\ -\int_{u=0}^t
e^{-u \tilde F}\ \frac{d\tilde F}{d\s } \ e^{-(t-u)\tilde F}\, du,
$$
and so
$$
\align &      \  \frac{\partial \tilde F^{-s-1} }{\partial \s
}\biggl|_{\s =0} \ =\  \frac{1}{\Gamma(s+1)} \int_0^\infty t^{s}
\frac{\partial\,e^{-t\tilde F}}{\partial \s }\biggl|_{\s =0}\, dt
\ =\ - \frac{1}{\Gamma(s+1)} \int_0^\infty  \int_0^t t^{s}
\  e^{-u F  }\  F^\prime  \ e^{-(t-u) F  }\, dudt \\
&\qquad\qquad \ =\  -\frac{1 }{\Gamma(s+1)} \ \int_{0}^\infty
\int_{0}^\infty\  (u+v)^{s}
  \ e^{-u  F }\    F^\prime   \  e^{-v F }\, du dv.
\endalign
$$
Hence
$$
\multline  \frac{d^2 \Cal Z(s)}{d \s ^2}
\\ =\ \frac{-\Gamma(s+1)}{\Gamma(s-n/2)}\trace F^\pprime F^{-s-1}
 \ +\  \frac{1 }{\Gamma(s-n/2)} \trace  \int_{u=0}^\infty
\int_{v=0}^\infty\  (u+v)^{s}
 F^\prime  \ e^{-u  F }\    F^\prime   \  e^{-v F }\, du dv .
\endmultline
$$
\medskip
\medskip

In the case when $F$ has non trivial smoothly varying kernel, we have
  $\Pi\, F= 0= F\,\Pi$, so
 we get the  relations
 $$
 F\,\Pi^\prime\ =\ -F^\prime\, \Pi,\qquad\qquad\qquad \Pi^\prime\, F\ =\ -\Pi\, F^\prime.
 $$
 They imply in particular that
 $$
 F\,\Pi^\prime\, F\ =\ 0,\qquad\qquad\qquad\qquad \Pi \,F^\prime \Pi\ =\ 0.
 $$
 Since $\Pi^\prime=\Pi^\prime\Pi+\Pi\Pi^\prime$, we also have
 $$
 \Pi\Pi^\prime\Pi\ =\ 0,
 $$
 so
 $$
 (F+\Pi) \,\Pi^\prime \,(F+\Pi)\ =\ F\,\Pi^\prime\, \Pi\ +\ \Pi \,\Pi^\prime
 \,F\ =\ - F^\prime\,\Pi\ -\ \Pi\, F^\prime.
 $$
Hence
$$
\Pi^\prime\ =\ -(F+\Pi)^{-1}\, F^\prime\, \Pi\ -\
\Pi\,F^\prime\,(F+\Pi)^{-1}\ =\ -F^{-1}\, F^\prime\, \Pi\ -\
\Pi\,F^\prime\,F^{-1}.
$$
We define
$$
\Psi_{s}(\mu)
\ =\ \cases (1-\mu^{-s-1})/(\mu-1)\ &\ \mu\neq 0,1\\ s+1 \ &\ \mu=1\\ 0 \ &\ \mu=0\endcases
$$
Then when $\mu> 0$,
$$
\int_0^\infty\int_{0}^\infty\  (u+v)^{s}\, e^{-u  }\, e^{-v\mu }\, du dv
\ =\ \Gamma(s+1)\Psi_s(\mu).
$$
Starting from the formula
$$
\frac{d \Cal Z(s)}{d \s }\ =\
\frac{-\Gamma(s+1)}{\Gamma(s-n/2)} \trace F^\prime (F+\Pi)^{-s-1}
$$
and following the previous argument gives
$$
\multline \frac{d^2 \Cal Z(s)}{d \s ^2}
\ =\ \frac{-\Gamma(s+1)}{\Gamma(s-n/2)}\trace F^\pprime (F+\Pi)^{-s-1}\\
 \ +\  \frac{1 }{\Gamma(s-n/2)} \trace  \int_{0}^\infty
\int_{0}^\infty\  (u+v)^{s}
 F^\prime  \ e^{-u  (F+\Pi) }\    (F+\Pi)^\prime   \  e^{-v (F+\Pi) }\, du dv .
\endmultline
\tag 3.5
$$
Furthermore,
$$
\multline
 \trace  \int_{0}^\infty
\int_{0}^\infty\  (u+v)^{s}
 F^\prime \ e^{-u  (F+\Pi) }\   \Pi^\prime\    \  e^{-v (F+\Pi) }\, du dv\\
=\ -2\trace  \int_{0}^\infty
\int_{0}^\infty\  (u+v)^{s}
 F^\prime \ e^{-u  (F+\Pi) }\   \Pi\ F^\prime\  F^{-1}  \  e^{-v (F+\Pi) }\, du dv\\
=\ -2\trace\left(
 F^\prime\   \Pi\ F^\prime\  F^{-1}  \ \int_{0}^\infty
\int_{0}^\infty\  (u+v)^{s} e^{-u}  e^{-v (F+\Pi) }\, du dv\right)\\
=\ -2\Gamma(s+1)\,\trace \ \Pi\ F^\prime\  F^{-1}  \ \Psi_{s}(F) F^\prime,
 \endmultline
$$
and since
$$
e^{-u(F+\Pi)} \ =\ e^{-uF}-\Pi\ +\ e^{-u}\Pi,
$$
we get
$$
\multline
  F^\prime \ e^{-u  (F+\Pi) }\   F^\prime \  e^{-v (F+\Pi) }
\ =\  F^\prime \ (e^{-u  F }-\Pi)\   F^\prime \ ( e^{-v F }-\Pi)\\ +\
 e^{-u}\ \ F^\prime \Pi\ F^\prime \  e^{-v F }
\ +\  e^{-v}   F^\prime e^{-uF} F^\prime \Pi,
\endmultline
$$
and
$$
\multline
\trace  \int_{0}^\infty
\int_{0}^\infty\  (u+v)^{s}
 F^\prime  \ e^{-u  (F+\Pi) }\    F^\prime   \  e^{-v (F+\Pi) }\, du dv\\
=\   \trace  \int_{0}^\infty
\int_{0}^\infty\  (u+v)^{s}
 F^\prime  \ (e^{-u  F }-\Pi)\    F^\prime   \  (e^{-v F }-\Pi)\, du dv\\
 \qquad\ +\ 2\Gamma(s+1)\,\trace \ \Pi F^\prime\,  \Psi_{s}(F)\, F^\prime.
 \endmultline
  $$
Moreover
$$
\trace F^\pprime \Pi\ =\  \trace (F^\prime \Pi)^\prime\ -\ \trace F^\prime \Pi^\prime
\ =\ -\trace F^\prime\Pi^\prime
\ =\ 2\trace \Pi F^\prime F^{-1}F^\prime,
$$
 and
  $$
2\trace \ \Pi F^\prime\,  \Psi_{s}(F)\, F^\prime
\ -\ 2\trace \ \Pi\ F^\prime\  F^{-1}  \ \Psi_{s}(F) F^\prime
\ -\ 2\trace \Pi F^\prime F^{-1}F^\prime\ =\ -2\trace \Pi F^\prime F^{-s-2} F^\prime.
  $$
Hence (3.5) becomes
$$
\multline \frac{d^2 \Cal Z(s)}{d \s ^2}
\ =\ \frac{-\Gamma(s+1)}{\Gamma(s-n/2)}\trace F^\pprime F^{-s-1}
\ -\  \frac{2\Gamma(s+1)}{\Gamma(s-n/2)}\trace \Pi F^\prime F^{-s-2} F^\prime\\
 \ +\  \frac{1 }{\Gamma(s-n/2)} \trace  \int_{0}^\infty
\int_{0}^\infty\  (u+v)^{s}
 F^\prime  \ (e^{-u  F }-\Pi)\    F^\prime   \ ( e^{-v F}-\Pi)\, du dv .
\endmultline
$$
We will reorganize this formula to combine terms which contain poles.  We have
$$
\multline
\trace  \iint_{u+v<1}\  (u+v)^{s}
 F^\prime  \ (e^{-u  F }-\Pi)\    F^\prime   \ ( e^{-v F}-\Pi)\, du dv\ = \\
 \trace  \iint_{u+v<1}\  (u+v)^{s}
 F^\prime  \ e^{-u  F }\    F^\prime    e^{-v F}\, du dv
 \ -\ 2  \trace  \iint_{u+v<1}\  (u+v)^{s}
 F^\prime  \Pi\    F^\prime     (e^{-v F}-\Pi)\, du dv    \\
 =\ \trace  \iint_{u+v<1} (u+v)^{s}
 F^\prime   e^{-u  F }\,    F^\prime     e^{-v F}\, du dv
 \ +\ 2  \trace
 F^\prime   \Pi\,    F^\prime  F^{-1} \left( \int_{0}^1
  t^{s} e^{-t F}\, dt\, -\, \frac{I}{s+1}\right).
 \endmultline
 $$
 So defining $L(s)$ as in (3.4), we obtain Lemma 3.1.

\medskip
\medskip
\medskip


\noindent{\bf 4. Reducing to the heat kernel at small times.}
\medskip

Our aim in this paper is to express the Hessian of $Z(s)$ in the
form given in (1.4). In this Section we reduce to a local problem
involving the heat kernel at small times. Most of the steps in
this reduction are quite routine.

We start by expressing $F$ defined in (1.2) in local coordinates.
Let $S^2_+(\Bbb R^n)$ be the cone of positive definite elements of
$S^2(\Bbb R^n)$.  There exist functions $a$, $b$, $c$, such that
in the natural trivialization of $E$ associated to arbitrary local
coordinates on $M$,
$$
 F \ =\  I\, V^{2/n} \,\sum_{ij} g^{ij}\partial_i
\partial_j \ +\ \sum_i  a^i \partial_i \ +\
b\ +\ \sum_{i,j,k,\ell} \left(\partial_k
\partial_\ell g^{ij} \right)c_{ij}^{k\ell} ,
\tag 4.1
$$
where $I$ is the identity $N\times N$ matrix,  and
$$
\align
 & a :(0,\infty)\times S^2_+(\Bbb R^n)\times T^{0,3}(\Bbb R^n)\ \to\
  T^{1,0}(\Bbb R^n)\otimes T^{1,1}(\Bbb R^N),\\
& b:(0,\infty) \times S^2_+(\Bbb R^n)\times T^{0,3}(\Bbb R^n)\
\to\
 \hphantom{T^{1,0}(\Bbb R^n)\otimes} T^{1,1}(\Bbb R^N),\\
 &  c:(0,\infty)\times S^2_+(\Bbb R^n)\hphantom{\times T^{0,3}(\Bbb R^n)}\ \ \,\to
 \ T^{2,2}(\Bbb R^n)\otimes T^{1,1}(\Bbb R^N)\\
 \endalign
 $$
are smooth functions which are independent of the metric and the
choice of local coordinates,
 so  writing $\partial
g$ for the tensor $(\partial_i g_{jk})\in T^{0,3}(\Bbb R^n)$ at each
point of the coordinate chart,
$$
\align
  \text{for}\qquad 1\leq i\leq n, \qquad\qquad & a^i=a^i(V,g,\partial g)
  \ \in\ T^{1,1}(\Bbb R^N),\\
 & b=b(V,g,\partial g)\ \in\ T^{1,1}(\Bbb R^N),\\
\text{for}\qquad 1\leq i,j,k,\ell\leq n, \qquad\qquad &
c_{ij}^{k\ell}=c_{ij}^{k\ell}(V,g)\ \in\ T^{1,1}(\Bbb R^N).
 \endalign
$$

In the previous Section we computed the
second variation of $Z(s)$ when the varying metric is given by
$$
\tilde g(\s)\ =\ g\ +\ \s h.
$$
Differentiating (4.1) once, we obtain
$$
F^\prime\ =\  \sum\Sb i,j\\{|\al|+|\beta|\leq 2}\endSb \,
(\partial^\al h_{ij})\,A^{ij}_{\al\beta} \,\partial^\beta,
\qquad\qquad\qquad
 A^{ij}_{\al\beta}  =A^{ij}_{\al\beta}( V,\partial^\gamma
 g)\ \in\ T^{1,1}(\Bbb R^N),\tag 4.2
$$
where  $\gamma$ ranges over multiindices with
$|\beta|+|\gamma|\leq 2$. Differentiating (4.1) twice we obtain
$$
\align &  F^\pprime\ =\ \sum\Sb i,j,k,\ell \\
|\al|+|\beta|\leq 2,\\ |\gamma|\leq2
 \endSb
(\partial^\al h_{k\ell})(\partial^\beta h_{ij})\,B_{\al\beta
\gamma }^{ij k\ell} \,\partial^\gamma, \qquad\qquad\qquad
B_{\al\beta \gamma}^{ ij k\ell}=B_{\al\beta \gamma}^{ ij
k\ell}(V,\partial^\delta g)\ \in\ T^{1,1}(\Bbb R^N).
\tag 4.3\\
\endalign
$$

The Hessian of $Z(s)$ with respect to the flat connection on the
space of metrics is just given by
$$
\Hess Z(s)(h,h)\ =\ \frac{d^2}{d\s^2} Z(s).
$$
We will show that if we can prove Theorem 1 in this case, then the
result also follows when we use the intrinsic connection on the
space of metrics. We also reduce (1.4) to the case when  $\llangle
g,h\rrangle_g =0$. We will then express the second variation given
in Lemma 3.1, in terms of $h$, and analyze the terms which arise.

\definition{Definition}
Write
 $\pi_j:M\times M\to M$ for the projection onto the $j^{th}$
factor. Fix a metric $g$ on $M$. For a pseudodifferential operator
$Q:C^\infty(E)\to C^\infty(E)$ of order less than $-n$, define
$K(Q,x,y)$ to be the continuous  section  of $\pi_1^*E\otimes
\pi_2^* E^*$ on $M\times M$ such that $K(Q,x,y)\,dV(y)$ is the
integral kernel of $Q$.
 If $Q(s)$ is an analytic family of polyhomogeneous pseudodifferential operators then
 $K(Q(s),x,y)$ may be  analytically continued from the $s$-region where
 $Q(s)$ has order less than $-n$, and we write $K(Q(s),x,y)$ to denote
 this analytic continuation.  When $x\neq y$, the section $K(Q(s),x,y)$
 is entire in $s$ and $K(Q(s),x,y)dV(y)$  gives the Schwartz
 kernel of $Q(s)$ away from the diagonal. At the diagonal,
 $K(Q(s),x,x)$ is a meromrphic in $s$.
 To compute  the poles and the trivial zeros of $K(Q(s),x,x)$, we refer the reader to  [Ok1]
 Sections 2 and 3.
\enddefinition

 From (3.3) and (4.2),
$$
\multline
 D  \Cal Z(s) (h)\ =\ \trace F^\prime L(s) \ =\
\int_M \trace K(F^\prime L(s) , x,x)\, dV(x)\\
=\  \sum\Sb i,j\\{|\al|+|\beta|\leq 2}\endSb \int (\partial^\al
h_{ij}(x))\,K(\,A^{ij}_{\al\beta}
 \,\partial^\beta L(s) , x,x)\, dV(x).
\endmultline
$$
By introducing a partition of unity and integrating by parts, we can
express this in the form
$$
  \llangle h\,,\, X(s)\rrangle_g.
  \tag 4.4
$$
where  $X(s)\in C^\infty(S^2M)$.  The expression for $L(s)$ in
(3.4) was chosen precisely so that (3.3) holds and
$$
K(A^{ij}_{\al\gamma} \partial^\al L(s) , x,x)
$$
is an entire function of $s$.

Now we reduce Theorem 1 to the case when the following two
conditions hold:

(a). We take the flat connection on $\Cal M $

(b). $\llangle h,g\rrangle_g =0$.

To make this reduction, assume that (1.4) holds when (a) and (b)
are satisfied.  We will first remove condition (b).  For general
$h\in C^\infty( S^2(TM))$, set
 $$
 \tilde{\tilde g}(\s,\tau)\ =\  e^\tau (g\ +\ \s h),
 $$
 so
 $$
 \partial_\s\tilde{\tilde g}(0,0)\ =\ h,
\qquad\qquad\qquad \partial_\tau\tilde{\tilde g}(0,0)\ =\ g.
 $$
 Then $\Cal Z(\tilde{\tilde g})$ is independent of $\tau$, so from (2.1)
 we see that
 $$
\Hess_\flat \Cal Z(s)( h,g)\ =\
-DZ(s)\left(\partial_\s\partial_\tau \tilde g\right)\ =\ -\ D\Cal
Z(s)(h).
$$
Define $P:C^\infty(S^2M)\to C^\infty(S^2M)$ by
$$
Ph=h-c_h g, \qquad\qquad \text{where}\qquad\qquad c_h \  =\ \frac1{nV}  \llangle h,g\rrangle_g,
$$
Then defining $X(s)$ as in (4.4),
$$
\align & \Hess_\flat \Cal Z(s)(h,h)\ =\ \Hess_\flat \Cal Z(s)(Ph, Ph)
\, +\, 2c_h\Hess_\flat \Cal Z(s)(h, g)\, +\, c_h^2\Hess_\flat \Cal Z(s)(g,g) \\
& =\ \llangle Ph , T_s Ph \rrangle_g\ -\ 2c_h D\Cal Z(s)(h)\\
&=\ \llangle Ph, T_s P h\rrangle_g\ -\
\tfrac2{nV}\llangle h,g\rrangle_g \llangle h, X(s)\rrangle_g\\
&=\ \llangle h, PT_s P h\rrangle_g\ -\
\Llangle\ h\ ,\ \tfrac2{nV}\llangle g,h\rrangle\ X(s) \Rrangle_g
\endalign
$$
 The operator corresponding to the first term is $PT_sP$ which
is a $\Psi$DO equal to $T_s$ up to the
smoothing operator
$T_s-PT_s P$.
The operator corresponding to the second term is $h\to \tfrac2{nV} \llangle g,h\rrangle X(s)$
which is smoothing and analytic in $s$.  This removes condition (b).

By (2.2), for the intrinsic connection on $\Cal M$ we have
$$
\Hess_\intrinsic \Cal Z(s)(h,h)\ =\ \Hess_\flat \Cal Z(s)(h,h)\ +\ D\Cal Z(s)(w)
$$
where
$$
w_{ij}\ =\ -\ \sum_{k,\ell} h_{ik}g^{k\ell}   h_{\ell j} \ +\ \frac12 \langle h,g\rangle_g
\,h_{ij}- \frac14\langle h,h\rangle_g\, g_{ij}
$$
Then raising indices with respect to the metric $g$ in the usual way,
$$
 D\Cal Z (s)(w)\ =\ -\int \sum_{i,j,k,\ell} h_{ik} g^{k\ell}  h_{ \ell j}X^{ij}(s)\, dV \ +\
\frac12\int_M \langle h,g\rangle_g \,\langle h, X(s)\rangle_g\,dV\
-\ \frac14 \int_M\langle h,h \rangle_g \,\langle g,X(s)\rangle_g\, dV
$$
The operator corresponding to this form is a degree zero operator
which is analytic in $s$.

This reduces Theorem 1 to the case when we take the
flat connection on $\Cal M $ and assume that $h$ is trace free
relative to $g$.
In particular we now have $\tilde
V^\prime(0)=0$.  The rest of this section is devoted to analyzing the
 terms of the Hessian of $\Cal Z(s)$ in Lemma 3.1, on the right hand side of (3.1).

With $L(s)$ defined in (3.4), from (4.3) we have
$$
\trace\,F^\pprime\,L(s)
 \ =\  \int_M  \sum\Sb i,j,k,\ell \\
|\al|+|\beta|\leq 2,\\ |\gamma|\leq2
 \endSb(\partial^\al h_{k\ell})(\partial^\beta
h_{ij})\,
 \trace K(B_{\al\beta \gamma }^{ij k\ell} \partial^\gamma\,L(s),x,x)\, dV(x).
 \tag 4.5
$$
The operator corresponding to this form has degree $2$ and is entire in $s$.

To deal with the remaining three terms in (3.1) we will need the
following standard  estimates on the asymptotics of the heat
kernel $K(e^{-tF},x,y)$. Write $d$ for the geodesic distance
between points $x$ and $y$, and write $\nabla$ for the connection
on $E$ or $E^*$ associated to $g$.  We take the norms $|\,\cdot\,
|_g$ associated to $g$ on the fibers  $T^{p,q}_{(x,y)}(M\times
M)$. There exists $\ep>0$ such that  for all $j,k,\ell$,
$$
\align
& \sup_{d(x, y)>d_0 }
\left| \partial_t^j \nabla_x^k \nabla_y^\ell K(e^{-tF},x,y)\right|
 \ \leq\ C e^{-\ep /t},\qquad\qquad\text{ when }t<1, \tag 4.6 \\
& \sup_{ x,y\in M}
\left| \partial_t^j \nabla_x^k  \nabla_y^\ell K(e^{-tF}-\Pi,x,y)\right|
 \ \leq\ C e^{-\ep t},\qquad\qquad \text{ when }t>1, \tag 4.7
\endalign
$$
where the constant $C$ depends on $j,k,\ell$.
Moreover, if we write
$\pi_i:[0,1]\times M\times M\to M$, $i=1,2$ for the projection
onto the $i$th  $M$ factor, and $d(x,y)$ for the distance from $x$ to $y$,
there exists
$ b\in C^\infty(
\pi_1^*E\otimes \pi_2^* E^*)$ and $d_0>0$, $\ep>0$ such that
$$
K(e^{-tF},x,y)\ =\ \frac1{(4\pi t)^{n/2}V} e^{-d^2/(4tV^{2/n})}
b(t,x,y), \qquad\qquad \text{ when }t<1, d<d_0. \tag 4.8
$$
Moreover, $b(0,x,x)$ is the identity matrix  in $E_x\otimes
E^*_x$.  See for example [Ch].

By (4.7), we see that
$$
\int_1^\infty t^s K(e^{-tF}-\Pi,x,y)\, dt
$$
is a smooth kernel which is entire in $s$, and so there exists a
smoothing operator $P(s)$ which is entire in $s$ such that
$$
\trace\left( F^\prime \Pi F^\prime F^{-1}
 \int_{1}^\infty\  t^{s}\, (e^{-t F}-\Pi)\, dt \right)
 \ =\ \llangle h, P(s)h\rrangle.\tag 4.9
 $$
Now consider the term from (3.1):
$$
\align
& \trace  \iint_{u+v>1}\  (u+v)^{s}
 F^\prime  \ (e^{-u  F }-\Pi)\    F^\prime   \ ( e^{-v F}-\Pi)\, du dv\\
&=\ \trace \int_1^\infty\int_1^\infty\qquad\qquad \quad\qquad\qquad"
\qquad \qquad\qquad\,dudv\\
&\ +\ \trace \int_{u=0}^{1/2}\int_{v=1-u}^\infty\qquad\qquad \quad\qquad"
\qquad\qquad \qquad\qquad\,dudv\\
&\ +\ \trace \int_{v=0}^{1/2}\int_{u=1-v}^\infty\qquad\qquad \quad\qquad"
\qquad\qquad \qquad\qquad\,dudv\\
&=\ \text{(I)}\ +\ \text{(II)}\ +\ \text{(III)}.\tag 4.10
 \endalign
 $$
In any of these terms we can pass the trace inside the integral.
Indeed, writing $\|\cdot\|$ for the operator norm,
$\|\cdot\|_1$ for the trace class norm and
and $\|\cdot\|_2$ for the Hilbert-Schmidt norm of an operator, when $u>1/2$,
$$
\multline
\|F^\prime  \ (e^{-u  F }-\Pi)\    F^\prime   \ ( e^{-v F}-\Pi)\|_1
\ \leq\   \|F^\prime  \ (e^{-u F/2 }-\Pi)\|_2  \| (e^{-u F/2 }-\Pi)
\ F^\prime \|_2 \|  \ ( e^{-v F}-\Pi)\|\\
 \leq\ Ce^{-\ep u},
 \endmultline
$$
where the last inequality follows from (4.7). This shows for
example that
$$
\trace \int_{v=0}^{1/2}\int_{u=1-v}^\infty (u+v)^{s}
 \|F^\prime  \ (e^{-u  F }-\Pi)\    F^\prime   \ ( e^{-v F}-\Pi)\|_1\, du dv \ <\ \infty.
 $$
 so that the trace can be passed through the integral in (III).

Now we localize by writing $1=\sum_{\chi\in \Lambda} \chi$ where
$\Lambda$ is a finite partition of unity on $M$ such that for each
point $x\in M$, all the functions $\chi\in\Lambda$ which contain
$x$ in their supports are supported in a single  contractible
coordinate chart, $\Omega$ over which $E$ is trivial, and  on
which $\partial ^\al g_{ij}$, and $\partial^\al g^{ij}$ are
uniformly bounded for each $\al$. Hence we express
$$
\multline \trace  F^\prime  \ (e^{-u  F }-\Pi)\    F^\prime   \ (
e^{-v
F}-\Pi) \\ =\ \sum\Sb \chi, \tilde\chi\in \Lambda,\ i,j,k,\ell\\
|\al|+|\beta|\leq 2,\ |\gamma|+|\delta|\leq 2\endSb \trace\ \chi
\,A_{\al \beta}^{ij}\ (
\partial^\al h_{ij})\partial^\beta \ \left(e^{-u F }-\Pi\right) \
\tilde\chi \, A_{\gamma \delta}^{k\ell}\ (
\partial^\gamma h_{k\ell})\partial^\delta  \ \left( e^{-v F
}-\Pi\right).
\endmultline\tag 4.11
$$
We can analyze each term in the above sum separately. Fix
$\chi,\tilde\chi, i,j,k,\ell,\al,\beta,\gamma,\delta$ and set
$$
\tilde E(u,v,x,y)\ =\ \chi(x)\tilde\chi(y)\ \trace\ A_{\al
\beta}^{ij}(x)\ (\partial^\beta_x K(e^{-u  F }-\Pi,x,y)) \
 A_{\gamma \delta}^{k\ell}(y)\  (\partial^\delta_y K(  e^{-v
F }-\Pi,y,x)).\tag 4.12
$$
The trace on the right is the trace on $N\times N$ matrices. Then
$$
\multline \trace\ \chi \,A_{\al \beta}^{ij}\ (
\partial^\al h_{ij})\partial^\beta \ \left(e^{-u F }-\Pi\right) \
\tilde\chi \, A_{\gamma \delta}^{k\ell}\ (
\partial^\gamma h_{k\ell})\partial^\delta  \ \left( e^{-v F
}-\Pi\right)
\\ =\
 \iint_{\Supp(\chi) \times\Supp(\tilde\chi)}
   ( \partial^\al h_{ij}(x)) \tilde E(u,v,x,y) ( \partial^\gamma h_{k\ell}(y))\,
   \sqrt{g(x)} \sqrt{g(y)}\,dxdy.
   \endmultline
   \tag 4.13
$$
 Now
from (4.7), for $u>1$ and $v>1$ and each $\al$, $\beta$,
$$
|\partial_x^\al \partial_y^\beta \tilde E(u,v,x,y)|\ \leq\
Ce^{-\ep(u+v)}.
$$
Hence
$$
 K_s^{(I)}(x,y)\ :=\  \frac{1 }{\Gamma(s-n/2)}  \int_{1}^\infty\int_{1}^\infty
 (u+v)^s \tilde E(u,v,x,y)\,dudv
 \tag 4.14
$$
is smooth in $(s,x,y)$ and entire in $s$ and (I) is a finite sum of terms of the form
$$
   \iint_{\Supp(A) \times\Supp(B)}
   ( \partial^\al h_{ij}(x)) K_s^{(I)}(x,y)
   ( \partial^\gamma h_{k\ell}(y))\,\sqrt{g(x)} \sqrt{g(y)}\, dxdy.
   \tag 4.15
$$
Integrating by parts to take the derivatives from $h$,
we find that the kernel of this bilinear form is
smooth and entire in $s$.
Similarly, write (II) formally as a sum of terms of the form
$$
\iint_{\Omega\times\Omega}  ( \partial^\al h_{ij}(x))\
\left(\chi(x)\ \trace K_{s}^{(II)}(x,y)\ \tilde\chi (y)\ \right) (
\partial^\gamma h_{k\ell}(y))\, \sqrt{g(x)} \sqrt{g(y)} \, dxdy,
\tag 4.16
$$
where $K_s^{(II)}(x,y)$ is given by
$$
K_s^{(II)}(x,y)\ =\ \int_{u=0}^{1/2} \int_{v=1-u}^\infty (u+v)^s
A^{ij}_{\al\beta}(x) (
\partial^\beta_x K(e^{-u  F },x,y)) \ A_{\gamma\delta}^{k\ell}(y)\ (\partial^\delta_y K(  e^{-v
F },y,x))\, dudv.
$$
In doing this we are switching the $dxdy$ integration with the
$dudv$ integration, which must properly interpreted and justified,
since  the  integral in (4.16)  may be divergent  due to
$K_s^{(II)}(x,y)$ being singular on the diagonal.  The way to
proceed is to work in the class of pseudodifferential operators of
order $|\beta|$.  Define
$$
J_s^{(II)}(u,x,y)\ :=\ A_{\gamma\delta}^{k\ell}(y)\
\int_{v=1-u}^\infty (u+v)^s (\partial^\delta_y K( e^{-v F },y,x))
\, dv\ A^{ij}_{\al\beta}(x)\,
$$
so
$$
K_s^{(II)}(x,y)\ =\   \int_{u=0}^{1/2} (
\partial^\beta_x K(e^{-u  F },x,y))\
\ J_s^{(II)}(u,x,y)\, du. \tag4.17
$$
From (4.7), $J_s^{(II)}(u,x,y)$  is smooth for $(s,u,x,y)\in \Bbb
C\times [0,\tfrac12]\times \Omega\times\Omega$ and analytic in
$s$.  We notice that the integrand on the right hand side of
(4.17) is the integral kernel of an operator on $C^\infty(\Omega)$
whose order $|\beta|$ symbol norms are bounded uniformly in the
parameter $u$. This follows easily from the elementary lemmas
below, and shows that the integral in (4.17) converges to the
Schwartz kernel of an order $|\beta|$ pseudodifferential operator.

\proclaim{Lemma 4.1} The operators $e^{-uF}$ are psedodifferential
operators of order zero, with symbol norms uniformly bounded in
the parameter $u$.
\endproclaim
    This implies that
for $\beta$ fixed, the operators $\partial^\beta_x e^{-u  F }$
are psedodifferential operators of order $|\beta|$ with uniformly
bounded symbol norms in the parameter $u$. \proclaim{Lemma 4.2} If
$\Cal K(Q,x,y)$ is the Schwartz kernel of a $\Psi$DO, $Q$, on $M$
of order $m$, and $J(x,y)$ is smooth.  Then the product $\Cal
K(Q,x,y)J(x,y)$ is the Schwartz kernel
 of a $\Psi$DO, which we denote by  $Q_J$, which also has order $m$.
   If $Q$ depends smoothly on a parameter $p\in\Bbb R^d$ then so does $Q_J$.
 \endproclaim
We can conclude that the expression in (II) is equal to (4.16),
and $K_s^{(II)}(x,y)$  is  the Schwartz kernel of a $\Psi$DO of
order $|\beta|$. The question then arises whether this operator is
polyhomogeneous and how one obtains its symbol expansion. For a
function $\psi\in C^\infty([0,1])$, integration by parts gives
$$
\multline
\int_{u=0}^{1/2} e^{-u\la} \psi(u)\, du
\ =\ \la^{-1}\psi(0)\ +\ \la^{-2} \psi^\prime(0)
\ + \ \dots\ +\ \la^{-M}\psi^{(M-1)}(0)
\\ -\ \left( \la^{-1}\psi(1)\ +\ \la^{-2} \psi^\prime(1)
\ + \ \dots\ +\ \la^{-M}\psi^{(M-1)}(1)\right) e^{-\la/2} \ +\
\la^{-M}\int_{u=0}^{1/2}  e^{-u\la} \psi^{(M)}(u)\, du,
\endmultline
$$
we see that (4.17) equals
$$
\multline
\sum_{m=0}^{M-1} ( \partial^\beta_x K(F^{-m-1} ,x,y))\
\partial_u^m|_{u=0} J_s^{(II)}(u,x,y)
\ -\ \sum_{m=0}^{M-1} ( \partial^\beta_x K(F^{-m-1}e^{-F/2} ,x,y))\
\partial_u^m|_{u=1} J_s^{(II)}(u,x,y)\\
\ +\ \int_{u=0}^{1/2} ( \partial^\beta_x K(F^{-M}e^{-u  F },x,y)) \
\partial_u^M J_s^{(II)}(u,x,y)\, du.
\endmultline
\tag 4.18
$$
The term in the second sum
$$
- ( \partial^\beta_x K(F^{-m-1}e^{-F/2} ,x,y))\ \partial_u^m|_{u=1}
J_s^{(II)}(u,x,y)
$$
is smooth in $(s,x,y)$ and analytic in $s$. We need the following
refinement of Lemma 4.2: \proclaim{Lemma 4.3} In fact if the
$\Psi$DO, $Q$, in Lemma 4.2 is polyhomogeneous of degree $\mu$,
then so is $Q_J$, and  if $J=J_s$ depends
 analytically on a parameter $s$, then $Q_{J_s}$ is an analytic
 family of pseudodifferential operators.
\endproclaim
With this refinement, we see that the term in (4.18)
$$
 ( \partial^\beta_x K(F^{-m-1} ,x,y))\
\partial_u^m|_{u=0} J_s^{(II)}(u,x,y)\, dV(y)
$$
is  the Schwartz kernel of a polyhomogeneous $\Psi$DO operator of
degree $|\beta|-2m-2$.  The last term in (4.18) is the Schwartz
kernel of a pseudodifferential operator of order $|\beta|-2M-2$.
We conclude that  $K_s^{(II)}(x,y)$  is  the Schwartz kernel of a
polyhomogeneous $\Psi$DO of order $|\beta|-2$, which is analytic
in $s$. The lemmas are consequences of the standard theory. Indeed
for Lemma 4.1, following [Gi] one can express the symbol of the
heat operator in terms of the symbol of the resolvant via a
contour integral. However, the resolvant $(\la-F)^{-1}$ is a
pseudodifferential operator, and the zeroth order symbol norms of
$\la(\la-F)^{-1}$ are uniformly bounded in $\la$. Lemmas 4.2 and
4.3 follow from the relationship between amplitudes and symbols,
see [Tr].

We return to our analysis of (4.10). Term (III) is dealt in the
same fashion as Term (II) with the roles of $u$ and $v$
interchanged. Now we want to understand the final term in (3.1)
which is given in (1.10):
$$
\frac{1 }{\Gamma(s-n/2)} \trace  \iint_{u+v<1}\  (u+v)^{s}
 F^\prime  \ e^{-u  F }\    F^\prime   \  e^{-v F}\, du dv.
 \tag 4.19
$$
From now on we will refer to the steps listed in Table 1.1.

\noindent{\bf Step 14}: Using estimates based on (4.8), we can
take the trace in (4.8) inside the integral because
$$
\|F^\prime  \ e^{-u  F }\    F^\prime   \  e^{-v F}\|_1
\ \leq\  \|F^\prime  \ e^{-u  F/2 }\|_2 \|e^{-uF/2}\    F^\prime \|_2 \|
 e^{-v F}\|\ \leq\ C u^{-n-2}.
$$
Since this is also true with $u$ and $v$ interchanged, we have
$$
\|F^\prime  \ e^{-u  F }\    F^\prime   \  e^{-v F}\|_1\ \leq\ C(u+v)^{-n-2},
$$
and so for $\Re s$ sufficiently large,
$$
 \iint_{u+v<1}\  (u+v)^{s}
\| F^\prime  \ e^{-u  F }\    F^\prime   \  e^{-v F}\|_1\, du dv\ <\ \infty.
 $$
{\bf Steps 7, 3}: Taking a partition of unity
$1=\sum_{\chi\in\Lambda}\chi$ on $M$ as for (4.11), we see that
$$
 \trace  F^\prime  \ e^{-u  F }\    F^\prime   \  e^{-v
F} \ =\ \sum\Sb \chi, \tilde\chi\in \Lambda,\ i,j,k,\ell\\
|\al|+|\beta|\leq 2,\ |\gamma|+|\delta|\leq 2\endSb \trace\ \chi
\,A_{\al \beta}^{ij}\ (
\partial^\al h_{ij})\partial^\beta \ e^{-u F } \
\tilde\chi \, A_{\gamma \delta}^{k\ell}\ (
\partial^\gamma h_{k\ell})\partial^\delta  \  e^{-v F}.
\tag 4.20
$$
{\bf Step 1}: Fixing $\chi,\tilde\chi\in\Lambda,
i,j,k,\ell,\al,\beta,\gamma,\delta$, we will consider each term in
(4.20). Set
$$
E(u,v,x,y)\ =\ \chi(x)\tilde\chi(y)\ \trace\ A_{\al
\beta}^{ij}(x)\ (\partial^\beta_x K(e^{-u  F },x,y)) \
 A_{\gamma \delta}^{k\ell}(y)\  (\partial^\delta_y K(  e^{-v
F },y,x)).\tag 4.21
$$
{\bf Step 10}: Then
$$
\multline \trace\ \chi \,A_{\al \beta}^{ij}\ (
\partial^\al h_{ij})\partial^\beta \ e^{-u F } \
\tilde\chi \, A_{\gamma \delta}^{k\ell}\ (
\partial^\gamma h_{k\ell})\partial^\delta  \  e^{-v F
}
\\ =\
 \iint_{\Supp(\chi) \times\Supp(\tilde\chi)}
   ( \partial^\al h_{ij}(x)) E(u,v,x,y) ( \partial^\gamma h_{k\ell}(y))\,
   \sqrt{g(x)} \sqrt{g(y)}\,dxdy.
   \endmultline
   \tag 4.22
$$
{\bf Steps 8, 14}: Hence formally, (4.19) is a sum of terms of the
form
$$
  \iint_{\Supp(\chi) \times\Supp(\tilde\chi)}
   ( \partial^\al h_{ij}(x)) K_s(x,y) ( \partial^\gamma h_{k\ell}(y))\,
   \sqrt{g(x)} \sqrt{g(y)}\, dxdy,
   \tag 4.23
$$
where
$$
 K_s(x,y)\ :=\  \frac{1 }{\Gamma(s-n/2)}  \iint_{u+v<1}
 (u+v)^s E(u,v,x,y)\,dudv.
 \tag 4.24
$$
When the supports of $\chi$ and $\tilde\chi$ are disjoint we can
use (4.6) and (4.8) to see that (4.19) is a bilinear form on
$C^\infty(S^2M)$ with a smooth kernel which is analytic in $s$.
Hence we have reduced to the case when $\chi$ and $\tilde \chi$
are supported on a coordinate chart $\Omega$ whose diameter is
less than the injectivity radius of $M$.
\medskip
\medskip
\medskip


\noindent{\bf 5. The symbol of $U_s$.}
\medskip

 Here in
Section 5, we  describe  how to compute the symbol expansion of
$U_s$, and we carry this out explicitly for $u_s$ thus  obtaining
Theorem 2.  To do this we expand on the summary  in Section 1.
Section 5 can be read independently from Section 6, which contains
 the rigorous proof of Theorem 1, and  shows that
the algorithm we lay down here is correct. We will see in
particular in Section 6, that (4.19) is a sum of terms of the form
(4.23). Moreover, we see in Lemmas 6.1 and 6.5 that $K_s$ defined
in (4.24) has symbol
$$
\tilde U_s(x,\xi)\ +\ \tilde V_s(x,\xi),
$$
where setting
$$
 m=|\beta|+|\delta|,
 $$
$\tilde U_s$ is  polyhomogeneous of degree $n-2s-4+m$,
 and $\tilde V_s$ polyhomogneous of degree $m-2$.

\noindent{\bf Step 3}: We are interested in obtaining the symbol
expansion of $\tilde U_s$ at the point $x\in M$, and so we work in
coordinates on a small coordinate chart $\Omega$ containing $x$.

\noindent{\bf Steps 11, 12, 13}: To obtain the symbol expansion of
$\tilde U_s$, we make an expansion of $E(u,v,x,y)$. Using (4.8),
$$
  \partial_{x}^\beta K(e^{-u  F }, x,y)
\ =\ u^{-n/2} e^{-d^2(x,y)/4u} b_\beta(u, x,y)\tag 5.1
$$
where $u^{|\beta|}b_\beta(u, x,y)$ is a smooth $N\times N$ valued
function on $[0,1]\times\Omega^2$ whose Taylor expansion  at $u=0$
gives an expansion of the form
$$
 b_\beta(u,x,y)\ \sim\ \sum_{\ell=-|\beta|}^\infty  \ \,u^\ell\
b_{\beta\ell}(x,y), \qquad\qquad b_{\beta\ell}(x,y)\ =\
O(\,d(x,y)^{\max\{-2\ell-|\beta|,0\}}\,).
$$
Here, $b_{\beta\ell}(x,y)$ is smooth on $\Omega^2$. Hence
considering (4.21),
$$
 E(u,v,x,y)\ =\
(uv)^{-n/2} e^{-d^2(x,y)(1/4u+1/4v)} J(u,v,x,y),\tag 5.2
$$
where
$$
u^{|\beta|}v^{|\delta|}J(u,v,x,y)
$$
 is smooth on $[0,1]^2\times\Omega^2$.
and $J(u,v,x,y)$ has an expansion of the form
$$
J(u,v,x,y)\ \sim\
\sum_{a=-|\beta|}^\infty\sum_{b=-|\delta|}^\infty u^a v^b
J_{ab}(x,y), \qquad\qquad J_{ab}(x,y)\ =\
O(\,d(x,y)^{\max\{-2(a+b)-(|\beta|+|\delta|),0\}}\,).\tag 5.3
$$
Here, $J_{ab}(x,y)$ is smooth on $\Omega^2$.

 \noindent{\bf Step
4}: Fix $x$ and write $w$ for  geodesic normal coordinates of $y$
centered at $x$ . We make a Taylor expansion of $J_{ab}(x,y)$
about $w=0$ to write
$$
J(u,v,x,y)\ \sim\ \sum_{\ell=-m}^\infty J_\ell(u,v,x,y), \tag 5.4
$$
where for constants $J_{ab\mu}(x)$,
$$
J_\ell(u,v,x,y)\ =\ \sum_{2a+2b+|\mu|=\ell}  J_{ab\mu} u^a v^b
w^\mu.\tag 5.5
$$
The sum in (5.5) is clearly finite.
 Now set
$$
\align & E_{\ell}(u,v,x,y)\ =\ (uv)^{-n/2}
e^{-d^2(x,y)(1/4u+1/4v)}
J_\ell(u,v,x,y),  \tag 5.6\\
& K_{s,\ell}(x,y)\ =\ \frac1{\Gamma(s-n/2)}\iint_{u+v<1} (u+v)^s
 E_{\ell}(u,v,x,y)\, dudv.
\endalign
$$
{\bf Step 15}: The symbol of $E_\ell(u,v,x,\xi)$ has the form
$$
U_{s,\ell}(x,\xi)\ +\ V_{s,\ell}(x,\xi),
$$
where $U_{s,\ell}(x,\xi)$ is homogeneous of degree $n-2s-4-\ell$,
and $V_s(x,\xi)$ is polyhomogeneous of degree $m-2$. To see this,
the metric-modified symbol (see 2.4) of   $E_\ell(u,v,x,y)$ is
obtained by taking the Fourier transform in $y$,
$$
 \s(E_\ell,u,v,x,\xi)\ =\ \int e^{i w\cdot \xi}
E_\ell(u,v,x,y)\, dw. \tag 5.7
$$
From the definition of $E_\ell$,
$$
\align
 \s & (E_\ell,u,v,x,\xi)\ =\ (4\pi)^{n/2} (u+v)^{-n/2}
\sum_{2a+2b+|\mu|=\ell} J_{ab\mu}(x) u^a v^b (-i\partial_\xi)^\mu
e^{-|\xi|^2uv/(u+v)} \\
&  =\  (u+v)^{-n/2}e^{-|\xi|^2uv/(u+v)}\sum\Sb 2a+2b+|\mu|=\ell\\
\nu\leq\mu \endSb c_{\mu\nu} J_{ab\mu}(x) (u+v)^{-(|\mu|+|\nu|)/2}
(uv)^{(|\mu|+|\nu|)/2} u^a v^b\xi^\nu
\tag 5.8  \\
& \\
 &=\ (u+v)^{-n/2-\ell}e^{-|\xi|^2uv/(u+v)}\sum\Sb |\nu|\leq 2m+\ell\\ 2a+2b=3\ell+|\nu|
\endSb \tilde J_{ab\nu}(x) u^a v^b \xi^\nu,
\endalign
$$
where the constants $c_{\mu\nu}$ and functions $\tilde
J_{ab\nu}(x)$ are defined appropriately by collecting terms. We
will prove in Lemma 6.3 that $\tilde J_{ab\nu}(x)$ vanishes if $a$
or $b$ is negative, so
$$
\s(E_\ell,u,v,x,\xi)\ =\ O((u+v)^{-n/2-\ell})
$$
as $u,v\to 0$.

\noindent{\bf Step 14}: Hence the symbol of $K_{s,\ell}(x,\xi)$ in
(4.24)  equals
$$
\frac1{\Gamma(s-n/2)}\sum\Sb |\nu|\leq \ell\\ 2a+2b=3\ell+|\nu|
\endSb \tilde J_{a b\nu}(x)\xi^\nu
\iint_{u+v<1}(u+v)^{s-n/2-\ell}u^a v^b e^{-|\xi|^2 uv/(u+v)}\,
dudv, \tag 5.9
$$
which converges absolutely for $\Re s $ large.

\noindent{\bf Step 16}: By making the change of variables
$$
T=(u+v)|\xi|^2,\qquad\qquad\qquad \tau=\frac u{u+v},
$$
we see in Lemmas 6.4 and 6.5 (also (1.19)) that this  integral
extends analytically and to a sum of symbols of degrees $m-2$ and
$n-2s-4-\ell$.  The symbol of degree $n-2s-4-\ell$ is obtained by
extending the integral on the right of (5.9) to
$(u,v)\in[0,\infty)^2$, and so
$$
\multline U_{s,\ell}(x,\xi)\ =\ \sum\Sb |\nu|\leq\ell\\
2a+2b=3\ell+|\nu|
\endSb \tilde J_{a b\nu}(x)\ \xi^\nu
\ |\xi|^{n-2s-4-|\nu|-\ell} \\
\times\qquad
\frac{\Gamma(s-n/2-\ell+a+b+2)\Gamma(n/2-s+\ell-b-1)\Gamma(n/2-s+\ell-a-1)}
{\Gamma(s-n/2)\Gamma(n-2s+2\ell-a-b-2)}.
\endmultline
\tag 5.10
$$
One might hope to complete the analysis of term (4.19) by studying
the properties of  the remainder term in the Taylor expansion
(5.4). However  although the order of $U_{s,\ell}$ decreases  as
$\ell$ increases, the order of $V_{s,\ell}$  does not. We will
need to decompose the kernel $E(u,v,x,y)$ in additional ways which
pick out the homogeneous terms in the symbol of $\tilde V_s$. This
is carried out in Section 6.  We emphasize that (5.10) gives the
term of degree $n-2s-4-\ell $ in the symbol  of $K_s(x,y)$. From
(4.23), $u_s(x,\xi)$ is the symbol of
$$
\sum\Sb \chi, \tilde\chi\in \Lambda,\ i,j,k,\ell\\
|\al|+|\beta|\leq 2,\ |\gamma|+|\delta|\leq 2\endSb
(-1)^{|\al|+|\gamma|}\ (\sqrt{g(x)}\sqrt{g(y)})^{-1}
\,\partial_x^\al\partial_y^\gamma\left(\sqrt{g(x)}\sqrt{g(y)}
K_s(x,y)\right), \tag 5.11
$$
which can be computed from (5.10) using the symbolic calculus. It
is also easy to go from the metric-modified symbol to the true
symbol using the symbolic calculus if one so desires.  See [Tr]
for details on the symbolic calculus. One needs to  be a little
careful because  we have been working in normal coordinates based
at $x$, and in order to apply the usual symbolic calculus formulas
to compose operators, one has use the same coordinate system for
the points $x$ and $y$. The necessary change of variables formulas
can also be found in [Tr].
\medskip
\medskip

\demo{Proof of Theorem 2}  We want to compute $u_s(x,\xi)$ which
is the degree $n-2s$ term in the symbol of (5.11). (As always, in
saying this we must assume $s$ is not a half-integer greater than
$n/2-1$, since at these values the symbol of $\tilde V(s)$ also
contributes to the degree-$(n-2s)$ term.) Now  $U_{s,-m}(x,\xi)$
is the degree $n-2s-4+m$ term in the symbol of $K_s(x,y)$, which
is the principal symbol when $\Re s$ is small.  Hence from (5.11),
$$
  u_s(x,\xi)\ =\
\sum\Sb \chi, \tilde\chi\in \Lambda,\ i,j,k,\ell\\
|\al|+|\beta|\leq 2,\ |\gamma|+|\delta|\leq 2\endSb
(-1)^{|\al|}(i\xi)^{\al+\gamma} U_{s,-m}(x,\xi). \tag 5.12
$$
We will evaluate 5.12 at a fixed point $x=x_0$, and from now on we
express $x$ and $y$ in  geodesic normal coordinates centered at
$x_0$.  In particular, we assume that the derivatives in (4.21)
are computed in these coordinates.

\proclaim{Lemma 5.1} Let $C(s)$ be defined as in (1.7). Define
$$
G\ =\  (4\pi)^{-n/2}(u+v)^{-n/2}\ \  e^{-|\xi|^2uv/(u+v)},\tag
5.13
$$
define the polynomial $P(\partial^\al,t,w)$ in the variable $w\in
\Bbb R^n$ by
$$
P(\partial^\al,t,w)\ =\  e^{|w|^2/4t }\partial^\al_w e^{-|w|^2/4t
},\tag 5.14
$$
define
$$
\s(\partial^\al,\,\partial^\beta,\,\partial^\gamma,\,\partial^\delta,\
u,v,x,\xi) \ =\ (i\xi)^{\al+\gamma} \
P(\partial^\beta,u,-i\partial_\xi)\
P(\partial^\delta,v,-i\partial_\xi)G, \tag 5.15
$$
and define
$$
u_s(\partial^\al,\,\partial^\beta,\,\partial^\gamma,\,\partial^\delta,\
x,\xi)\ =\ \frac{1 }{C(s)\Gamma(s-n/2)}
\int_{0}^\infty\int_{0}^\infty
 (u+v)^s   \s(\partial^\al,\,\partial^\beta,\,\partial^\gamma,\,\partial^\delta,\
u,v,x,\xi)  \,dudv.\tag 5.16
$$
Then
$$
u_s(x_0,\xi)\ =\ C(s)\sum\Sb i,j,k,\ell\\
|\al|+|\beta|= 2,\ |\gamma|+|\delta|= 2\endSb
u_s(\partial^\al,\,\partial^\beta,\,\partial^\gamma,\,\partial^\delta,\,x_0,\xi)\trace
A_{\al\beta}^{ij}(x_0) A_{\gamma\delta}^{k\ell}(x_0).
 \tag 5.17
$$
\endproclaim

\demo{Proof of Lemma 5.1} Recall that each term $K_s(x,y)$ in
(4.24) depends on the parameters $\chi,\tilde\chi\in\Lambda,
i,j,k,\ell,\al,\beta,\gamma,\delta$, and is given by
$$
 K_s(x,y)\ :=\  \frac{1 }{\Gamma(s-n/2)}  \int_{0}^\infty\int_{0}^\infty
 (u+v)^s E(u,v,x,y)\,dudv.
$$
where
$$
E(u,v,x,y)\ =\ \chi(x)\tilde\chi(y)\ \trace\ A_{\al
\beta}^{ij}(x)\ (\partial^\beta_x K(e^{-u  F },x,y)) \
 A_{\gamma \delta}^{k\ell}(y)\  (\partial^\delta_y K(  e^{-v
F },y,x)).
$$
The degree $n-2s-4+m$ term in the symbol of  $K_s$ is
$U_{s,-m}(x,\xi)$, where
$$
 U_{s,-m}(x,\xi) \ =\ \frac1{\Gamma(s-n/2)}
\int_0^\infty\int_0^\infty (u+v)^{s}\s(E_{-m},u,v,x,\xi)    \,
dudv,\tag 5.18
$$
and  $E_{-m}$ which gives the lowest order terms of $E$ is defined
in (5.6).  We will see that writing $w$ for the geodesic normal
coordinates of $y$ about $x=x_0$, we have
$$
\multline
 E_{-m}(u,v,x_0,y)\\ =\ \chi(x_0)\tilde\chi(x_0)\
(4\pi)^{n/2}(uv)^{-n/2}(-1)^{|\beta|} (\partial^\beta_w
e^{-|w|^2/4u })
 \ (\partial^\delta_w   e^{-|w|^2/4v})\ \trace A_{\al\beta}^{ij}(x_0)
  A_{\gamma\delta}^{k\ell}(x_0).
  \endmultline
 \tag 5.19
$$
To see this, from (5.1) we must compute the ``lowest order" terms
in the functions
$$
 e^{d^2(x,y)/4u} \partial_x^\beta \left( e^{-d^2(x,y)/4u}
b(u,x,y)\right)\biggl|_{x=x_0}, \qquad\qquad e^{d^2(x_0,y)/4v}
\partial_y^\delta \left( e^{-d^2(x_0,y)/4v} b(v,x_0,y)\right),
$$
about $u=v=0$ and $y=x_0$, where the orders of the terms $u^a
w^\mu$ and $v^a w^\mu$ are  $2a+|\mu|$.  It is clear that to get
the lowest order terms, we can replace $b(u,x,y)$ and $b(v,x_0,y)$
by $b(0,x_0,x_0)$ which is the identity matrix. This reduces the
problem to computing the lowest order terms in the expansions of
the functions
$$
\align &  e^{d^2(x,y)/4u} \partial_x^\beta
e^{-d^2(x,y)/4u}\bigl|_{x=x_0},\tag 5.20 \\
& e^{d^2(x_0,y)/4v}
\partial_y^\delta e^{-d^2(x_0,y)/4v}.
\tag 5.21
\endalign
$$
Now since the derivatives are expressed in geodesic coordinates
about $x_0$,  (5.21) is precisely
$$
\ e^{|w|^2/4v}\
\partial_w^\delta e^{-|w|^2/4v},
$$
which is a linear combination of terms $v^a w^\mu$ with order
$2a+|\mu|=-|\delta|$.  On the other hand (5.20) may involve terms
coming from the curvature of the metric. However, these are all of
higher order than $-|\beta|$, as one can check by expressing $ x$
and $y$ in normal coordinates around $x_0$, and expanding
$$
d(x, y)^2\ =\ | x-y|^2\ +\ O\left(| x-y|^3\right).
$$
We thus find that the lowest order terms in (5.20) are given by
$$
  e^{|x-y|^2/4u}\ \partial_{x}^\beta e^{-|
x-y|^2/4u}\ \bigl|_{ x=x_0}\ =\  (-1)^{|\beta|}\ e^{|w|^2/4u}\
\partial_{w}^\beta e^{-|w|^2/4u} .
$$
From this, we get (5.19).  Since the Fourier transform of
$(4\pi)^{n/2}(uv)^{-n/2} e^{-|y|^2(1/4u+1/4v) }$ is $G$, taking
the Fourier transform of $E_{-m}$ gives
$$
\s(E_{-m},u,v,x_0,\xi)\ =\ \chi(x_0)\tilde\chi(x_0)\ \trace
A_{\al\beta}^{ij}(x_0) A_{\gamma\delta}^{k\ell}(x_0)\
(-1)^{|\beta|} P(\partial^\beta,u,-i\partial_\xi)\
P(\partial^\delta,v,-i\partial_\xi)G.
 \tag 5.22
$$
Plugging this into (5.12) and (5.18) and summing over $\chi,\tilde
\chi\in \Lambda$ yields Lemma 5.1.
\enddemo

To complete the proof of Theorem 2, it remains to compute
$u(\partial^\al,\,\partial^\beta,\,\partial^\delta,\,\partial^\gamma,\
u,v,x,\xi)$ for explicit values of $\al,\beta,\gamma,\delta$. Now
$$
\align
& P(\partial_{j},u,w) \ =\ -\frac{w_j}{2u}  ,\qquad\qquad\qquad\qquad\qquad\qquad
P(\partial_{j},u,-i\partial_{\xi}) \ =\ \frac{i}{2u}\partial_{\xi_j}.  \\
& P(\partial_{j}\partial_{k},u,w) \ =\  \left( \frac{w_jw_k }{4u^2}
\ -\ \frac{\delta_{jk}}{2u}\right)\qquad\qquad\qquad
P(\partial_{j}\partial_{k},u,-i\partial_\xi) \ =\  -\left( \frac{1 }{4u^2}\partial_{\xi_j}\partial_{\xi_k}
\ +\ \frac{\delta_{jk}}{2u}\right), \\
\endalign
$$
and
$$
\align
 &\partial_{\xi_j} G\ =\  - \frac{2uv}{u+v}\xi_j G, \\
 &\partial_{\xi_j}\partial_{\xi_k} G
 \ =\  \left( \frac{4u^2v^2}{(u+v)^2}\xi_j\xi_k \ -\
 \frac{2uv}{u+v}\delta_{jk}\right)G, \\
 & \partial_{\xi_j}\partial_{\xi_k} \partial_{\xi_p} G\ =\
 \left(-\frac{8u^3v^3}{(u+v)^3}\xi_j\xi_k\xi_p
 \ +\   \frac{4u^2v^2}{(u+v)^2}(\delta_{jk}\xi_p+\delta_{jp}\xi_k+\delta_{kp}\xi_j) \right)G.\\
 & \partial_{\xi_j}\partial_{\xi_k} \partial_{\xi_p}\partial_{\xi_q} G\ =\
 \biggl(\frac{16u^4v^4}{(u+v)^4}\xi_j\xi_k\xi_p\xi_q\ +\
 \frac{4u^2v^2}{(u+v)^2}(\delta_{jk}\delta_{pq}+\delta_{jp}\delta_{kq}+\delta_{kp}\delta_{jq})
 \\
 &\qquad\qquad\qquad\qquad\qquad\ -\   \frac{8u^3v^3}{(u+v)^3}(
 \delta_{jq}\xi_k\xi_p
 +\delta_{kq}\xi_j\xi_p
 +\delta_{pq}\xi_j\xi_k
 +\delta_{jk}\xi_p\xi_q
 +\delta_{jp}\xi_k\xi_q
 +\delta_{kp}\xi_j\xi_q
 )\biggr)G.\\
 \endalign
 $$
So for (a),
$$
\s(\partial_j\partial_k, \, I,\, \partial_p\partial_q,\,I, \ u,v,
x,\xi) \ =\  \xi_j\xi_k\xi_p\xi_q\, G.
$$
for (b),
$$
\multline \s(\partial_j, \, \partial_k,\,
\partial_p\partial_q,\,I,\ u,v, x,\xi) \ =\ -i\xi_j\xi_p\xi_q
P(\partial_k,u,-i\partial_\xi)G\ =\
 \frac{\xi_j\xi_p\xi_q}{2u}\partial_{\xi_k}   G \ =\   -\frac{v}{u+v}\xi_j\xi_k\xi_p\xi_q G.
\endmultline
$$
for (c),
$$
\multline \s(\partial_j, \, \partial_k,\, \partial_p,\,
\partial_q,\ u,v, x,\xi) \ =\ -\xi_j\xi_p
P(\partial_k,u,-i\partial_\xi)P(\partial_q,v,-i\partial_\xi)G\ =\
 \frac{\xi_j\xi_p}{4uv}\partial_{\xi_k} \partial_{\xi_q}   G \\
 =\    \left( \frac{uv}{(u+v)^2}\xi_j\xi_k\xi_p\xi_q \ -\
 \frac{ 1}{2(u+v)}\delta_{kq}\xi_j\xi_p\right)G.
\endmultline
$$
for (d),
$$
\multline \s(I, \, \partial_j\partial_k,\,
\partial_p\partial_q,\,I,\ u,v, x,\xi) \ =\ -\xi_p\xi_q
P(\partial_j\partial_k,u,-i\partial_\xi)G\ =\ \left(
\frac{\xi_p\xi_q }{4u^2}\partial_{\xi_j}\partial_{\xi_k}
\ +\ \frac{\xi_p\xi_q\delta_{jk}}{2u}\right)G \\
 =\  \frac{\xi_p\xi_q }{4u^2}\left( \frac{4u^2v^2}{(u+v)^2}\xi_j\xi_k \ -\
 \frac{2uv}{u+v}\delta_{jk}\right)G\ +\ \frac{\xi_p\xi_q\delta_{jk}}{2u}G
\ =\  \frac{v^2}{(u+v)^2}\xi_j\xi_k\xi_p\xi_q G\ +\
\frac1{2(u+v)}\delta_{jk}\xi_p\xi_q G.
\endmultline
$$
for (e),
$$
\multline \s(I, \, \partial_j\partial_k,\, \partial_p,\,
\partial_q,\ u,v, x,\xi) \ =\ i\xi_p
P(\partial_j\partial_k,u,-i\partial_\xi)P(\partial_q,v,-i\partial_\xi)
G
\\ =\ \left( \frac{\xi_p }{8u^2v}\partial_{\xi_j}\partial_{\xi_k}\partial_{\xi_q}
\ +\ \frac{\xi_p\delta_{jk}}{4uv}\partial_{\xi_q}\right)G
\\ =\ -\frac{uv^2}{(u+v)^3}\xi_j\xi_k\xi_p\xi_q G
 \ +\   \frac{v}{2(u+v)^2}(\delta_{jk}\xi_p\xi_q+\delta_{jq}\xi_k\xi_p+\delta_{kq}\xi_j\xi_p) G
\ -\  \frac{1}{2(u+v)}\delta_{jk}\xi_p \xi_qG\\
\endmultline
$$
and finally for (f),
$$
\align \s &(I,  \, \partial_j\partial_k, \,I,\, \partial_p
\partial_q,\ u,v, x,\xi) \ =\
P(\partial_j\partial_k,u,i\partial_\xi)P(\partial_p\partial_q,v,i\partial_\xi)
G
\\  &=\ \left( \frac{1}{4u^2}\partial_{\xi_j}\partial_{\xi_k}
\ +\ \frac{\delta_{jk}}{2u}\right)\left( \frac{1}{4v^2}\partial_{\xi_p}\partial_{\xi_q}
\ +\ \frac{\delta_{pq}}{2v}\right)G\\
 &=\ \left( \frac{1}{16u^2v^2}\partial_{\xi_j}\partial_{\xi_k}
\partial_{\xi_p}\partial_{\xi_q}
\ +\ \frac{\delta_{jk}}{8uv^2}\partial_{\xi_p}\partial_{\xi_q}
\ +\ \frac{\delta_{pq}}{8u^2v}\partial_{\xi_j}\partial_{\xi_k}\ +\ \frac{\delta_{jk}\delta_{pq}}{4uv}
\right)G\\
&=\ \biggl(\frac{u^2v^2}{(u+v)^4}\xi_j\xi_k\xi_p\xi_q\ +\
 \frac{1}{4(u+v)^2}(\delta_{jk}\delta_{pq}+\delta_{jp}\delta_{kq}+\delta_{kp}\delta_{jq}) \\
 &\qquad\qquad\qquad\qquad
 \ -\   \frac{uv}{2(u+v)^3}(\delta_{jq}\xi_k\xi_p
 +\delta_{kq}\xi_j\xi_p
 +\delta_{pq}\xi_j\xi_k
 +\delta_{jk}\xi_p\xi_q
 +\delta_{jp}\xi_k\xi_q
 +\delta_{kp}\xi_j\xi_q
 )\biggr)G\\
 &\ +\left( \frac{u}{2(u+v)^2}\delta_{jk}\xi_p\xi_q \, -\,
 \frac{1}{4(u+v)v}\delta_{jk}\delta_{pq}
 \, +\, \frac{v}{2(u+v)^2}\xi_j\xi_k\delta_{pq} \, -\,
 \frac{1}{4u(u+v)}\delta_{jk}\delta_{pq}\, +\, \frac1{4uv}\delta_{jk}\delta_{pq}\right)G\\
&=\ \biggl(\frac{u^2v^2}{(u+v)^4}\xi_j\xi_k\xi_p\xi_q\ +\
 \frac{1}{4(u+v)^2}(\delta_{jk}\delta_{pq}+\delta_{jp}\delta_{kq}+\delta_{kp}\delta_{jq})
 \ +\  \frac{u}{2(u+v)^2}\delta_{jk}\xi_p\xi_q\\
 &\ +\ \frac{v}{2(u+v)^2}\xi_j\xi_k\delta_{pq}
 \ -\   \frac{uv}{2(u+v)^3}(\delta_{jq}\xi_k\xi_p
 +\delta_{kq}\xi_j\xi_p
 +\delta_{pq}\xi_j\xi_k
 +\delta_{jk}\xi_p\xi_q
 +\delta_{jp}\xi_k\xi_q
 +\delta_{kp}\xi_j\xi_q)\biggr)G.\\
\endalign
$$
Now we just need to substitute these results for $\s$ into (5.16)
to obtain the expressions for $u_s$ given in Theorem 2 (a)-(f).
This is done more efficiently by looking at each type of term
which occurs in $\s$ and plugging it into (5.16).  Writing
$S=s-n/2$, we have
$$
\multline \frac1{\Gamma(s-n/2)} \int_0^\infty \int_0^\infty
(u+v)^{S-r}u^p v^q G\, dudv\ =
\\
  (4\pi)^{-n/2} |\xi|^{-2S+2r-2p-2q-4}
\frac{\Gamma(S-r+p+q+2)\Gamma(-S+r-p-1)\Gamma(-S+r-q-1)}{\Gamma(S)\Gamma
(-2S+2r-p-q-2)}\ = \\  \frac{S\cdots(S+1+p+q-r)\ \
(2S-1)\cdots(2S+2+p+q-2r)}{S\dots(S+1+p-r)\ \
S\cdots(S+1+q-r)}C(s) |\xi|^{-2S+2r-2p-2q-4}
\endmultline
$$
where $C(s)$ defined in (1.7), and  we assume $r\leq \min\{p,
q\}+2$. Hence we obtain
$u_s(\partial^\al,\,\partial^\beta,\,\partial^\delta,\,\partial^\gamma,\
x,\xi)$ from
$\s(\partial^\al,\,\partial^\beta,\,\partial^\delta,\,\partial^\gamma,\
u,v,x,\xi)$
 by substituting
$$
\align \frac{u^2 v^2}{(u+v)^4} G \ \rightarrow\
 (S^2+S)\, |\xi|^{n-2s-4},\qquad\qquad \qquad\quad\  \hfill &
  \frac{uv}{(u+v)^3} G \ \rightarrow\
 S\,  |\xi|^{n-2s-2},
   \\
\allowdisplaybreak \frac{u v^2}{(u+v)^3} G \ \rightarrow\
(2S^2+S-1) |\xi|^{n-2s-4},
 \ \qquad\qquad\quad \hfill &  \frac{v}{(u+v)^2}  G \ \rightarrow\
(2S-1)\, |\xi|^{n-2s-2},
  \\
\allowdisplaybreak
 \frac{uv}{(u+v)^2} G \
\rightarrow\
 (4S^2+2S-2) \,  |\xi|^{n-2s-4},
\qquad\qquad\  \ \hfill & \frac{1}{u+v} G \ \rightarrow\
 (4S-2)\, |\xi|^{n-2s-2},
  \\
 \allowdisplaybreak
\frac{v^2}{(u+v)^2}  G \ \rightarrow\
 (4S^2-2S) |\xi|^{n-2s-4},
 \qquad\ \  \qquad\qquad \hfill &  \frac{1}{(u+v)^2} G \ \rightarrow\
  |\xi|^{n-2s},
 \\
\allowdisplaybreak G \ \rightarrow\
 4(4S^2-1) |\xi|^{n-2s-4}
.\qquad\qquad\qquad \qquad\qquad\quad\ \hfill &  \\
\endalign
$$
\enddemo
\medskip
\medskip
\medskip


{\bf 6. The double decomposition.}
\medskip
\medskip

In this Section we complete the proof of Theorem 1. In order to
simplify the notation we will restrict  to the scalar case with
$V=1$.  It is straight forward  to obtain the vector case by
replacing scalars by matrices in the appropriate places.  Let
$\Omega$ be a coordinate chart of $M$ of small diameter which is
identified with a subset of $\Bbb R^n$, and make a change of
coordinates $(x,y)\to (x,w)$, sending $\Omega\times\Omega$ onto
$\Lambda\subset \Bbb R^{2n}$, where $(x,y)$ are product
coordinates on $\Omega\times\Omega$ and $w=w(x,y)$ are normal
coordinates of $y$ centered at $x$, which vary smoothly with $x$.
We define the function $\rho(x,w)$ to give the volume element on
$M$;
$$
\rho(x,w)\, dw\ =\ dV(y) =\ \sqrt{g(y)}\, dy.
$$
For a function $f$ on $\Omega\times\Omega$, we define $f^{\Cal N}$ on $\Lambda$ by
$$
f(x,y)\ =\ f^{\Cal N}(x,w(x,y)).
$$
 and we are
interested in understanding (4.19), which is a sum of terms of the
form
$$
 \frac{1 }{\Gamma(s-n/2)}  \iint_{u+v<1}
  \iint_{\Omega\times\Omega} (u+v)^s
   ( \partial^\al h_{ij}(x)) E(u,v,x,y) ( \partial^\gamma h_{k\ell}(y))\,\sqrt{g(x)} \sqrt{g(y)}\, dxdy\,dudv,
   \tag 6.1
$$
where
$$
E(u,v,x,y)\ =\ A(x)\,(\partial^\beta_x K(e^{-u  F },x,y)) \ B(y)\
(\partial^\delta_y K(  e^{-v F },y,x)). \tag 6.2
$$
Here, the smooth functions $A$ and $B$ are supported in $\Omega$.
Set
$$
m\ =\ |\beta|+ |\delta|. \tag 6.3
$$
{\bf Step 9}:   We denote the operator on $\Omega$ with kernel
$E(u,v,x,y)$ by $\Cal E_{u,v}$ and
 the metric-modified symbol of this operator (see (2.4)) by $\s(E,u,v,x,\xi)$, so
$$
\s(E,u,v,x,\xi)\ =\ \int_{\Bbb R^n} e^{iw\cdot\xi} E^{\Cal
N}(u,v,x,w)\, dw, \tag 6.4
$$
and for $f\in C^\infty(\Omega)$ compactly supported,
$$
(\Cal E_{u,v} f)(x)\ =\ \int_\Omega E(u,v,x,y) f(y)\sqrt{g(y)}\, dy
\ =\ \frac1{(2\pi)^n}\int_{\Bbb R^n}\int_\Omega \s(E,u,v,x,\xi) e^{-iw\cdot \xi} f^\Cal N(x,w)\rho(x,w)\, dw d\xi.
$$
{\bf Steps 14, 15}: We will see in Lemma 6.3, that
$$
|\s(E,u,v,x,\xi)|\ =\ O( (u+v)^{-n/2-m}),\qquad\qquad\qquad \text{ as }u,v\to0.
$$
Hence
$$
(\Cal E_{u,v} f)(x)\ =\  O( (u+v)^{-n/2-m}),\qquad\qquad\qquad \text{ as }u,v\to0,
$$
and when $\Re s>n/2+m-2$, we can interchange the $(u,v)$
integration in (6.1) with the $x$ integration and conclude that
(6.1) equals
$$
  \int_{\Omega}
    (\partial^\al h_{ij}) \  (\Cal W_s  \partial^\gamma h_{k\ell})\,\sqrt{g} dx.
    \tag 6.5
$$
where $\Cal W_s$ is the   operator
$$
 \frac{1 }{\Gamma(s-n/2)} \iint_{u+v<1} (u+v)^s \Cal E_{u,v}\,
dudv,
$$
which has symbol
$$
\s_s(x,\xi)\ :=\ \frac{1 }{\Gamma(s-n/2)} \iint_{u+v<1} (u+v)^s
\s(E,u,v,x,\xi)\, dudv. \tag 6.6
$$
 Our aim is to show the following Lemma. \proclaim{Lemma 6.1}
$\s_s(x,\xi)$ is entire in $s$ and has the form $ U_s+ V_s$ where
$ U_s$ has degree $n-2s-4+m$, $ V_s$ has degree $-2$, and $ U_s$
and $ V_s$ are meromorphic with simple poles in $n/2+\Bbb N^+$.
\endproclaim
This completes the proof of Theorem 1, since we can now  write
(6.5) in the form
$$
 \int_{\Omega}
     h_{ij} \  (\Cal T_s   h_{k\ell})\,\sqrt{g} dx.
$$
where $\Cal T_s$ is the sum of polyhomogeneous $\Psi$DOs of
degrees $n-2s$ and $4-m-2\leq 2$. We pause to outline the strategy
to prove Lemma 6.1.
 It is convenient to make the change of variables
$$
t\ =\ u+v,\qquad\qquad\qquad\qquad \tau=u/(u+v), \tag 6.7
$$
so
$$
u=t\tau,\qquad\qquad\qquad\qquad v=t(1-\tau).
$$
Write
$$
H(t,\tau,x,y)\ =\ E(t\tau, t(1-\tau),x,y), \tag 6.8
$$
so
$$
\s_s(x,\xi)\ =\ \frac{1 }{\Gamma(s-n/2)} \int_0^1\int_0^1 t^{s+1}
\s(H,t,\tau,x,\xi)\, d\tau dt. \tag 6.9
$$
By (4.8), we see that
$$
H(t,\tau,x,y)\ =\ t^{-n-m} (\tau(1-\tau))^{-n/2-m}
e^{-d^2(x,y)/(4t\tau(1-\tau))} P(t,\tau,x,y) \tag 6.10
$$
where
$$
P(t,\tau,x,y)\ \in\ C^\infty([0,1]^2\times \Omega^2). \tag 6.11
$$
{\bf Steps 12, 13}: We make a Taylor expansion of
$P(t,\tau,x,y)=P^{\Cal N}(t,\tau,x,w)$ in $w$. The $L^{\text{th}}$
order remainder term has order $|w|^L$, and in Lemma 6.2 we see
that it contributes a symbol of order $-L+C$ to $\s_s(x,\xi)$,
where $C=C(m,n,s)$  is independent of $L$. The consequence of this
is that we can ignore the Taylor remainder term and prove the
result when $P(u,v,x,y)$ is replaced by a Taylor polynomial.

\noindent{\bf Step 16}: Doing this and carrying out the $\tau$
integral in (6.9), we find in Lemma 6.4   that $\s_s(x,\xi)$ is a
finite sum of terms of the form
$$
\xi^\nu \int_{0}^1 t^{s+j} Q(t,|\xi|^2t,x)\, dt, \tag 6.12
$$
where $Q(t,T,x)$ is smooth on $[0,1]\times[0,\infty]\times \Omega$
(including  the point $T=+\infty$), and vanishes at $T=+\infty$.
We show in Lemma 6.5, that such a symbol is entire in $s$. Making
a Talyor expansion around $t=0$ and $T=\infty$ enables us to
obtain the polyhomogeneous expansions of $U_s$ and $ V_s$. Terms
of increasing power in $t$ give terms of decreasing order for
$U_s$, while terms of decreasing power in $T$ give terms of
decreasing order for $V_s$. We proceed to carry out this plan.

\demo{Proof of Lemma 6.1} \proclaim{Lemma 6.2} Suppose that $\s_s$
is given by (6.9) where $H$ is given by (6.10) and $P$ satisfies
(6.11). Fix $s\in\Bbb C$, and $k\in\Bbb N$, and suppose
  $$
 P^\Cal N(t,\tau,x,w)\ =\ O(|w|^L)
 $$
where
 $$
 L\ >\  k\ +\ (n-2\Re s-2)_++2m-2.
$$
Then $\s_s(x,\xi)$ is a symbol of order $-k$.
\endproclaim
This lemma is proved later on.
It is a standard fact in the theory of $\Psi$DOs that given a
a sequence $\s_j(x,\xi)$ of symbols on $\Omega$  with $s_j(x,\xi)$ homogeneous of degree
$\mu-j$, then there exists a symbol $\s(x,\xi)$ with symbol expansion
$$
\sum_{j=1}^\infty \s_j(x,\xi).
$$
Moreover, if each term $\s_j$ depends analytically on a parameter
$s$, then $\s$ can be chosen to depend analytically on $s$.
  This combined with Lemma 6.2 enables us to restrict attention
to the case when $P^\Cal N(u,v,x,w)$ is replaced by its Taylor polynomial in $w$ of degree $L$;
$$
P_L(t,\tau,x,w)\ :=\ \sum_{|\mu|\leq L} P_\mu(t,\tau,x) w^\mu.
$$
{\bf Steps 12, 13}: Set
$$
H_L(t,\tau,x,w)\ =\ t^{-n-m} (\tau(1-\tau))^{-n/2-m}
e^{-|w|^2/(4t\tau(1-\tau))} \sum_{|\mu|\leq
L}P_\mu(t,\tau,x)w^\mu. \tag6.13
$$
{\bf Step 15}: Then
$$
\align \s(H_L,t,\tau,x,\xi)\ &=\ \int_{\Bbb R^n} e^{i w\cdot \xi}
H_L(t,\tau,x,w)\, dw\\
 &=\ t^{-n-m} (\tau(1-\tau))^{-n/2-m}\sum_{|\mu|\leq L} P_\mu(t,\tau,x)\int_{\Bbb R^n} e^{i w\cdot \xi}
w^\mu e^{-|w|^2/(4t\tau(1-\tau))}  \, dw\\
&=\  (4\pi)^{n/2} t^{-n/2-m} (\tau(1-\tau))^{-m}
\sum_{|\mu|\leq L}P_\mu(t,\tau,x)  (-i\partial_\xi)^\mu
e^{-|\xi|^2 t\tau(1-\tau)} \\
 &=\  t^{-n/2-m}(\tau(1-\tau))^{-m}e^{-|\xi|^2 t\tau(1-\tau)}\sum_{|\mu|\leq L}
 P_\mu(t,\tau,x)\sum\Sb {|\nu|\leq |\mu|}\\
 |\mu+\nu|\ \text{even}\endSb
 c_{\mu\nu}(t\tau(1-\tau))^{|\mu+\nu|/2} \xi^\nu, \tag 6.14
\endalign
$$
where  $c_{\mu\nu}$ are constants. The problem with this expression is the presence of
the non-integrable factor
$$
(\tau(1-\tau))^{-m},
$$
which may not be neutralized by the factor
$$
(\tau(1-\tau))^{(|\mu+\nu|)/2}
$$
  when $|\nu|$ is small. However, the following lemma shows that this is not an issue after the terms
  are summed up.  We see from (6.14) that we can write $\s(H_L,t,\tau,x,\xi)$ in the form
  $$
\s(H_L,t,\tau,x,\xi)\ =\ t^{-n/2-m} e^{-|\xi|^2
t\tau(1-\tau)}\sum_{|\nu|\leq L} R_{L,\nu}(t,\tau,x)
  \xi^\nu.
$$
 \proclaim{Lemma 6.3} If $L\geq 3m$,  then
$$
R_{L,\nu}(t,\tau,x)
\ \in\  C^\infty([0,1]\times[0,1]\times\Omega).
$$
\endproclaim
This lemma is proved later on.  We now need to improve our
understanding of the dependence of $\s(H_L,t,\tau,x,\xi)$ on $t$.
Indeed,  Paying closer attention to what happens when we
differentiate the heat kernel, we see that  (6.13) can be improved
to
$$
 P_\mu(t,\tau,x)\ =\ O(t^{(m-|\mu|)_+/2}),
$$
and so from (6.14),
$$
 R_{L,\nu}(t,\tau,x)\ =\ O(t^{ (m+|\nu|)/2}).
$$
Hence writing
$$
\lceil r\rceil\ :=\ \inf \{ q\in \Bbb Z: q\geq r\},
$$
we see that
$$
\tilde R_{L,\nu}(t,\tau,x)\ :=\ t^{-\lceil (m+|\nu|)/2\rceil}R_{L,\nu}(t,\tau,x)\ \in\ C^\infty([0,1]^2\times\Omega).
$$
Now up to a term of order $-L+C$,
$$
\s_s(x,\xi)\ =\ \frac{1 }{\Gamma(s-n/2)} \sum_{|\nu|\leq L}\xi^\nu
\int_0^1\int_0^1 t^{s-n/2+\lceil(|\nu|-m)/2 \rceil+1} e^{-|\xi|^2
t\tau(1-\tau)}\tilde R_{L,\nu}(t,\tau,x) \, d\tau dt. \tag 6.15
$$
It is convenient to compactify the interval $[0,\infty)$ by including the point $+\infty$.
The differentiable structure at $+\infty$ is given by the coordinates
$\phi:(0,+\infty]\to [0,+\infty)$,
$$
\phi(T)\ = \ \cases 1/T  &T<+\infty\\ 0 & T=+\infty\endcases
$$
\proclaim{Lemma 6.4}  If
$$
R(t,\tau,x)
\ \in\  C^\infty([0,1]\times[0,1]\times\Omega).
$$
then
$$
Q(t,T,x)\ :=\ \int_0^1
e^{-T\tau(1-\tau)}R(t,\tau,x) \,
d\tau
$$
is smooth on $[0,1]\times[0,+\infty]\times \Omega$
and vanishes at  $T=+\infty$.
\endproclaim
This lemma is proved later on.  We need one more result to
complete the proof of Lemma 6.1.

\proclaim{Lemma 6.5 } Suppose $Q(t,T,x)$ is smooth on
$[0,1]\times[0,+\infty]\times\Omega$ and vanishes at $T=+\infty$.
Then
$$
\s(s,x,\xi)\ :=\   \int_{0}^1 t^{s-n/2} Q(t,|\xi|^2 t,x)\, dt
$$
is a symbol of order $-2$ on the set $\Re s>n/2-1$,
which extends analytically
to a meromorphic function for $s\in \Bbb C$, and
$$
\frac1{\Gamma(s-n/2+1)}\s(s,x,\xi)
$$
is smooth in $(s,x,\xi)$ and entire in $s$.
Moreover,
$$
\s(s,x,\xi)\ =\  U_s(x,\xi)\ +\  V_s(x,\xi),
$$
where
$$
\frac1{\sin(\pi s)} U_s(x,\xi),\qquad\qquad\qquad
\frac1{\sin(\pi s)} V_s(x,\xi),
$$
are polyhomogeneous of degrees $n-2s-2$ and $-2$ respectively,
smooth in $(s,x,\xi)$ and entire in $s$.  The symbol expansions
of $ U_s$ and $ V_s$ can be computed in terms of the Taylor series of $Q$.
Indeed, if
$$
\align
&Q(t,T,x)\ \sim\ \sum_{k=0}^\infty t^k q_k(T,x),\qquad t\to0,\\
&Q(t,T,x)\ \sim\ \sum_{\ell=0}^\infty T^{-\ell-1}q^\ell(t,x),\qquad T\to\infty
\endalign
$$
then as $|\xi|\to\infty$,
$$
\align
&U_s(x,\xi)\ \sim\ \sum_{k=0}^\infty u_{s,k}(x)|\xi|^{n-2s-2-2k},
\qquad\text{where}\qquad u_{s,k}(x)\ =\ \int_0^\infty T^{s-n/2+k}q_k(T,x)\, dT,\\
&V_s(x,\xi)\ \sim\ \sum_{k=0}^\infty v_{s,k}(x)|\xi|^{-2-2k},
\qquad\text{where}\qquad v_{s,k}(x)\ =\ \int_0^1 t^{s-n/2-k-1}q^k(t,x)\,dt.
\endalign
$$
\endproclaim
We apply Lemmas  6.4 and 6.5 to (6.15).  In order to apply Lemma
6.5 we replace $s$ by $s+\lceil (|\nu|-m)/2\rceil+1$ and deduce
that $\s_s(x,\xi)$ is analytic for $s>n/2+m/2-2$, and is a sum of
polyhomogeneous symbols of degrees $-2$ and $n-2s+m-4$ which are
meromorphic  in $s$ with poles in $n/2+\Bbb N^+$. Since $m$ is an
integer between $0$ and $4$,
 this completes the proof of Lemma 6.1. \qed\enddemo

The rest of this paper is devoted to proving Lemmas 6.2, 6.3, 6.4
and 6.5.
%
%
\demo{Proof of Lemma 6.2} We consider the kernel $K_s$ with symbol
$\s_s$.  This is given by
$$
\align
K_s(x,w)\ : &=\ \frac1{\Gamma(s-n/2)} \int_0^1\int_0^1 t^{s+1}H^\Cal N(t,\tau,x,w)\, dw\\
&=\ \frac1{\Gamma(s-n/2)}\int_0^1\int_0^1
t^{s+1-n-m} (\tau(1-\tau))^{-n/2-m}
e^{-|w|^2/(4t\tau(1-\tau))} P(t,\tau,x,w)  \,d\tau dt.\\
\endalign
$$
Now suppose $\kappa$ is a multiindex and $|\kappa|=k$.  Then
$$
\multline
\left|\partial_{w}^\kappa \partial_x^\nu K_s(x,w)\right|  \\
\  =\ \left|\frac1{\Gamma(s-n/2)}
\int_0^1\int_0^1    t^{s+1-n-m} (\tau(1-\tau))^{-n/2-m}
  \partial_w^\kappa\left( e^{-|w|^2/(4t\tau(1-\tau))} \partial_x^\nu P(t,\tau,x,w)\right) \,d\tau dt\right| \\
 \endmultline
$$
Using Leibnitz  rule for differentiation, we see that this
is bounded by a linear combination of terms of the following form where  $\ell$ is in the range $0\leq \ell\leq k$.
$$
  |w|^{L+2\ell-k} \int_0^1\int_0^1    t^{\Re s+1-n-m-\ell} (\tau(1-\tau))^{-n/2-m-\ell}
e^{-|w|^2/(4t\tau(1-\tau))}  \,d\tau dt,
$$
and setting $T=|w|^{-2}t$, this has the form
$$
 |w|^{L-k+2\Re s-2n-2m+4 }\int_0^{|w|^{-2}}   T^{\Re s+1-n-m-\ell} \int_0^1(\tau(1-\tau))^{-n/2-m-\ell}
e^{-1/(4T\tau(1-\tau))}  \,d\tau dT .
$$
Now
$$
\int_{0}^1  (\tau(1-\tau))^{-n/2-m-\ell}
e^{-1/(4T\tau(1-\tau)) } \, d\tau\ =\  \cases O(e^{-1/(2 T)}), \qquad & T\to 0 \\
 O(T^{n/2+m+\ell-1}) \qquad \qquad & T\to \infty, \endcases
$$
so we conclude that
$$
\left|\partial_{w}^\kappa \partial_x^\nu K_s(x,w)\right|\ \leq\
C(n,s)|w|^{L-k+2\Re s-2n-2m+4 }(1+|w|^{n-2\Re s-2}\log(c/|w|)),
$$
where $c$ is a constant larger than the diameter of $M$.
This implies that $\partial_{w}^\kappa \partial_x^\nu K_s(x,w)$ is integrable when
$$
L\ >\  k\ +\ (n-2\Re s-2)_++2m-2. \tag 6.16
$$
Now suppose that (6.16) holds for some integer $k\geq0$. Then
$$
\s_s(x,\xi)\ =\ \int_{\Bbb R^n} e^{i w\cdot \xi} K_s(x,w)\, dw,
$$
and
$$
\partial_x^\nu\partial_\xi^\mu \s_s(x,\xi)\  =\  \int_{\Bbb R^n} e^{i w\cdot \xi}
 (iw)^\mu \partial_x^\nu K_s(x,w)    dw.
$$
For fixed $\xi$,  there exists a $j$ such that
$$
| \xi_j|\ \geq \  \frac{|\xi|}{\sqrt{n}}.
$$
Then define  the multiindex $\kappa$ by
$$
\cases \kappa_\ell\ =\ 0  \quad  & \ell\neq j,\\
\kappa_j=k+|\mu|.\quad &\endcases,
$$
so
$$
|\xi^\kappa| \ \geq\ C\ |\xi|^{k+|\mu|},\qquad\qquad\text{where}\qquad C=C(n,k,\mu)>0.
$$
Since
$$
e^{iw\cdot \xi}=\frac{1}{(i\xi)^\kappa} (\partial_w^\kappa e^{i w\cdot \xi}),
$$
we can integrate
by parts to obtain
$$
\partial_x^\nu\partial_\xi^\mu \s_s(x,\xi)
 \  =\ \frac1{(-i\xi)^\kappa}\int_{\Bbb R^n}  e^{i w\cdot \xi}\
 \partial_w^\kappa  \partial_x^\nu \left((iw)^\mu  K_s(x,w) \right)\,d\tau dt dw.
$$
Replacing $L$ and $k$ in (6.16) by $L+|\mu|$ and $k+|\mu|$, and
choosing $R$ greater than the diameter of $M$,
$$
|\partial_x^\nu\partial_\xi^\mu \s_s(x,\xi)|\ \leq\
\frac{C(n,s,k,\mu)}{|\xi|^{k+|\mu|}}\int_{|w|\leq R} |w|^{L-k+2\Re s-2n-2m+4 }
\left(1+|w|^{n-2\Re s-2}\log(c/|w|)\right)  \, dw.
$$
The integral on the right is  finite  when (6.16) holds. This
completes the proof of  Lemma 6.2. \qed\enddemo


\demo{Proof of Lemma 6.3}  We will start by listing the relevant
definitions. Set
$$
\Lambda\ =\ \{(x,w):(x,y)\in\Omega^2\}\ \subset\ \Omega\times\Bbb R^n.
$$
We have
$$
 E(u,v,x,y)\ =\ A(x)\,(\partial^\beta_x K(e^{-u  F },x,y)) \ B(y)\ (\partial^\delta_y K(  e^{-v
F },y,x)),
$$
$$
m\ =\ |\beta|+|\delta|.
$$
Under the change of coordinates $(u,v)\to(t,\tau)$ the function $E$ transforms to $H$ which has the form,
$$
 H(t,\tau,x,y)\ =\ t^{-n-m} (\tau(1-\tau))^{-n/2-m}
e^{-d^2(x,y)/(4t\tau(1-\tau))} P(t,\tau,x,y) .
$$
where
$$
P(t,\tau,x,y)\ \in\ C^\infty([0,1]^2\times \Omega^2).
$$
Taking the $L^{th}$ order Taylor polynomial of $P^\Cal N(t,\tau,x,w)$ in $w$ we have
$$
 P_L(t,\tau,x,w)\ =\ \sum_{|\mu|\leq L} P_\mu(t,\tau,x) w^\mu.
 $$
 The part of $H^\Cal N$ corresponding to $P_L$ is $H_L$,
 $$
 H_L(t,\tau,x,w)\ =\ t^{-n-m} (\tau(1-\tau))^{-n/2-m}
e^{-|w|^2/(4t\tau(1-\tau))} \sum_{|\mu|\leq
L}P_\mu(t,\tau,x)w^\mu. \tag 6.17
$$
Taking the Fourier transform in $w$, we get the symbol $\s(H_L)$,
$$
\align
 \s(H_L,t,\tau,x,\xi)\ &=\  t^{-n/2-m}e^{-|\xi|^2 t\tau(1-\tau)}(\tau(1-\tau))^{-m}
 \sum_{|\mu|\leq L}P_\mu(t,\tau,x)\sum\Sb {|\nu|\leq |\mu|}\\ |\nu+\mu|\ \text{even}\endSb
 c_{\mu\nu}(t\tau(1-\tau))^{|\mu+\nu|/2} \xi^\nu \\
 &=\ t^{-n/2-m}e^{-|\xi|^2 t\tau(1-\tau)}\sum_{|\nu|\leq L}R_{L,\nu}(t,\tau,x)
  \xi^\nu.
  \endalign
  $$

Now we start the proof of Lemma 6.3 by showing that if $L>3m$ then
$$
R_{L,\nu}(t,\tau,x)\ \in\ C^\infty([0,1]\times(0,1]\times\Omega).
$$
This is equivalent to showing that
$$
t^{n/2+m}
e^{|\xi|^2 t\tau(1-\tau)}\s(H_L,t,\tau,x,\xi) \in\ C^\infty([0,1]\times(0,1]\times\Omega\times\Bbb R^n).
 \tag 6.18
$$
The key here is to write
$$
\multline
E^\Cal N(t,\tau,x,\xi)\\ =\ \sum_
{|\mu+\nu+\iota+\kappa|\leq m}C_{\mu,\nu,\iota,\kappa}\
 \partial_x^\mu\partial_w^\nu\bigl(\ S_{\mu,\nu,\iota,\kappa}(x,w)\ \
(\partial_x^\iota\partial_w^\kappa  K^\Cal N(e^{-t\tau  F },x,w)) \  \ K^\Cal N(  e^{-t(1-\tau)
F },w,x)\ \bigr),
\endmultline
\tag 6.19
$$
where
 $$
 S_{\mu,\nu,\iota,\kappa}(x,w)\ \in \ C^\infty([0,1]^2\times \Lambda),
 $$
 and the partial derivatives
are taken with respect to the $(x,w)$ coordinates. This is
accomplished by first changing variables from  $(x,y)$ to $(x,w)$,
and then using the product rule to take the derivatives off the
factor $K^\Cal N(  e^{-v F },w,x))$.  Now  from (4.8),
$$
K^\Cal N( e^{-uF}, x, w)\ =\ (4\pi u)^{-n/2} e^{-|w|^2/4u} b^\Cal N(u,x,w).
$$
Hence
$$
E^\Cal N(t,\tau,x,\xi)\ =\ \sum_{|\nu|\leq m}    t^{-n-m+|\nu|}\
\tau^{-n/2-m+|\nu|}\ (1-\tau)^{-n/2} \ \partial_w^\nu\bigl(\
e^{-|w|^2/(4t\tau(1-\tau))}\ X_\nu(t,\tau,x,w) \ \bigr), \tag 6.20
$$
where
$$
X_\nu\ \in\ C^\infty([0,1]^2\times \Lambda).
$$
Denote the $L^{th}$ order Taylor polynomial of $X$ in $w$ by
$$
X^L_\nu(t,\tau,x,w)\ =\ \sum_{|\mu|\leq L} X_{\nu,\mu}(t,\tau,x) w^\mu.
$$
Now  replacing $X_\nu$ in (6.20) by $X^L_\nu$ we form $H^L$,
$$
 H^L(t,\tau,x,w)\ :=\ \sum_{|\nu|\leq m}    t^{-n-m+|\nu|}\  \tau^{-n/2-m+|\nu|}\ (1-\tau)^{-n/2}
\ \partial_w^\nu\bigl(\
e^{-|w|^2/(4t\tau(1-\tau))}\ X_\nu^L(t,\tau,x,w) \ \bigr).
$$
Comparing this with $H_L$,
we find that only terms with degree close to $L$ contribute to the difference
$$
H^L(t,\tau,x,w)\ -\ H_L(t,\tau,x,w).
$$
Indeed, this difference has the form
$$
t^{-n-m}  (\tau(1-\tau))^{-n/2-m}
\ \ e^{-|w|^2/(4t\tau(1-\tau))}\ \sum_{L-m\leq |\th|\leq L+m} Y_\th(t,\tau,x) w^\th,
$$
where
$$
Y_\th(t,\tau,x)\ \in\ C^\infty([0,1]^2\times \Omega).
$$
Then computing the symbol by following (6.14),
$$
 \s(H^L,t,\tau,x,\xi)\ - \s(H_L,t,\tau,x,\xi)
 $$
 is a linear combination of terms of the  form
 $$
  t^{-n/2-m} (\tau(1-\tau))^{-m}  \sum_{L-m\leq |\th|\leq L+m} Y_\th(t,\tau,x)   (-i\partial_\xi)^\th
e^{-|\xi|^2 t\tau(1-\tau)}
$$
and this is a linear combination of terms of the form
$$
  t^{-n/2}(t\tau(1-\tau))^{-m+(|\th|+|\nu|)/2}e^{-|\xi|^2 t\tau(1-\tau)}
 Y_\th(t,\tau,x)  \xi^\nu,
$$
where
$$
L-m\leq |\th|\leq L+m,
\qquad\qquad\qquad
 (|\th|+|\nu|)/2\ \in\ \Bbb N.
$$
When $L\geq 3m$, we  conclude that $-m+(|\th|+|\nu|)/2\geq 0$ and so
$$
t^{n/2+m}e^{|\xi|^2 t\tau(1-\tau)}
( \s(H^L,t,\tau,x,\xi)\ - \s(H_L,t,\tau,x,\xi))\ \in\ C^\infty([0,1]^2\times \Omega\times\Bbb R^n).
$$
However, again following (6.14), we find that
$$
\s(H^L,t,\tau,x,\xi)
$$
is a linear combination of terms of the form
$$
  t^{-n/2-m}\  \tau^{-m} \  X_{\nu\mu}(t,\tau,x)\  (i\xi)^\nu
\  e^{-|\xi|^2 t\tau(1-\tau)},\qquad\qquad\qquad  X_{\nu\mu}\in C^\infty([0,1]^2\times\Omega).
$$
Hence
$$
t^{n/2+m}e^{|\xi|^2 t\tau(1-\tau)}
  \s(H^L,t,\tau,x,\xi) \ \in\ C^\infty([0,1]\times(0,1]\times \Omega\times\Bbb R^n).
$$
which proves (6.18). Interchanging the roles of $\tau$ and
$1-\tau$ by taking  the derivatives off the factor  $K^\Cal
N(e^{-t\tau  F },x,w))$  in (6.19), we prove similarly that
$$
t^{n/2+m}e^{|\xi|^2 t\tau(1-\tau)}
  \s(H_L,t,\tau,x,\xi) \ \in\ C^\infty([0,1]\times[0,1)\times \Omega\times\Bbb R^n).
$$
 This completes the proof of  Lemma 6.3.\qed
\enddemo


\demo{Proof of Lemma 6.4}  It is clear that since $R$ is smooth on
$[0,1]^2\times\Omega$,
$$
Q(t,T,x)\ =\ \int_0^1
e^{-T\tau(1-\tau)}R(t,\tau,x) \,
d\tau
$$
is smooth on $[0,1]\times[0,\infty)\times\Omega$. The  behavior  of $Q$ at $T=\infty$
is an application of the method of steepest descent.  Indeed,
$$
Q(t,T,x)\ =\ \int_{|\tau-1/2|<1/4}e^{-\tau(1-\tau)T}R(t,\tau,x) \,
d\tau\ +\ \int_{|\tau-1/2|>1/4}e^{-\tau(1-\tau)T}R(t,\tau,x) \,
d\tau.
$$
The first term on the right is smooth for  $(t,1/T,x)\in [0,1]\times[0,\infty)\times \Omega$ and
has zero Taylor expansion at $1/T=0$.  For the second term, make the change of variables
$\kappa=T\tau(1-\tau)$ to get
$$
\multline
Q(t,T,x)\ =\ \frac1T\ \int_{0}^{3T/16}  e^{-\kappa }\left(R(t\,,\ \tfrac12-\sqrt{\tfrac14-\kappa/T} \, ,\ x)
\ +\ R(t\,,\ \tfrac12+\sqrt{\tfrac14-\kappa/T} \, ,\ x)\right)\sqrt{\tfrac14-\kappa/T}\  d\kappa
\\ =\ \frac1T\ \int_{0}^{3T/16}  e^{-\kappa }f(t,\kappa/T,x)\, d\kappa
\endmultline
$$
where $f\in C^\infty([0,1]\times[0,3/16]\times \Omega)$. It is
straight forward to show that setting $\ep=1/T$, this is smooth
for $(t,\ep,x)\in [0,1]\times [0,\infty)\times\Omega$ and vanishes
at $\ep=0$. This completes the proof of Lemma 6.4.\qed\enddemo


\demo{Proof of Lemma 6.5}  We will need some additional notation.
We make  Taylor expansions of the terms $q_k(T,x)$ as
$T\to\infty$:
$$
q_k(T,x)\ \sim\ \sum_{\ell=0}^\infty q_k^\ell(x)T^{-\ell-1}.
$$
Then naturally the terms in these expansions also occur in the expansions of $q^\ell(t,x)$ as $t\to0$,
$$
q^\ell(t,x)\ \sim\ \sum_{k=0}^\infty q_k^\ell(x)t^k.
$$
We also have a Taylor expansion of $q_k(T,x)$ as $T\to 0$,
$$
q_k(T,x)\ \sim\ \sum_{m=0}^\infty q_{k,m}(x)T^m.
$$
Set
$$
Q_K(t,T,x)\ =\ \frac1{t^K}\left(Q(t,T,x)\ -\ \sum_{k=0}^{K-1}q_k(T,x)t^k\right),
$$
so $Q_k(t,T,x)$ is smooth up to $t=0$, and  as $T\to\infty$ we have the Taylor expansion
$$
Q_K(t,T,x)\ \sim\ \sum_{\ell=0}^\infty Q_K^\ell(t,x)T^{-\ell-1}
$$
where
$$
Q_{K}^\ell(t,x)\  =\  \frac1{t^K}\left(q^\ell(t,x)\ -\ \sum_{k=0}^{K-1}q_k^\ell(x)t^k\right).
$$
We will start by showing how to obtain the asymptotic expansion of
$\s(s,x,\xi)$ which holds in the $L^\infty$ norm.
The symbol estimates will be discussed afterwards.
Now $\s(s,x,\xi)$ converges and is analytic in $s$ for $\Re s>n/2-1$ and we obtain an analytic continuation by
using the Taylor series of $Q(t,T,x)$ at $t=T=0$.
Indeed, Taylor expanding $Q(t,T,x)$ at $t=0$,
$$
\align
&\s(s,x,\xi)\ =\ \int_{t=0}^1 t^{s-n/2+K}Q_{K}(t,t|\xi|^2,x) \, dt
\ +\ \sum_{k=0}^{K-1} \int_0^1 t^{s-n/2+k}q_k(|\xi|^2t,x)\, dt  \tag 6.21\\
\endalign
$$
The first integral converges for $\Re s>n/2-1-K$.  Changing variables,
$$
  \int_0^1 t^{s-n/2+k}q_k(|\xi|^2t,x)\, dt
\ =\ |\xi|^{n-2s-2-2k}\int_0^{|\xi|^2}T^{s-n/2+k}q_k(T,x)\, dT.
$$
Fixing $k<2L$, we will analyze
$$
\psi_k(s,x,\xi)\ :=\  \int_0^{|\xi|^2}T^{s-n/2+k}q_k(T,x)\, dT.
$$
This has an analytic continuation which is seen by making a Taylor expansion at
$T=0$.  Indeed,
$$
\multline
\psi_k(s,x,\xi)
\ =\ \int_0^{|\xi|^2}T^{s-n/2+K}\left( q_k(T,x)-\sum_{m=0}^{K-1} q_{k,m}(x)T^m\right) \, dT
\\ +\ \sum_{m=0}^{K-1}\frac{q_{k,m}(x)}{s-n/2+k+m+1}|\xi|^{2s-n+2k+2m+2}.
\endmultline
$$
The integral on the right hand side  converges and is analytic in $s$
for $\Re s>n/2-K-1$, and the other terms are
meromorphic with simple poles in the set $n/2-1-\Bbb N^+$. This completes the proof that
$\s(s,x,\xi)$ has a analytic continuation. Furthermore, set
$$
\phi_k(s,x,\xi)\ =\   \int_{|\xi|^2}^\infty T^{s-n/2+k}q_k(T,x)\, dT.
$$
This is convergent for $\Re s<n/2-k$, and it has an analytic continuation, which can be
seen by taking a
Taylor series as $T\to\infty$.  Indeed,
$$
\multline
\int_{|\xi|^2}^\infty T^{s-n/2+k}q_k(T,x)\, dT\ =\
\int_{|\xi|^2}^\infty T^{s-n/2+k}\left(q_k(T,x)-\sum_{\ell=0}^{L-1} q_k^\ell(x)T^{-\ell-1}\right)\, dT
\\ +\ \sum_{\ell=0}^{L-1}\frac{q_k^\ell(x)}{s-n/2+k-\ell}\ |\xi|^{2s-n+2k-2\ell}.
\endmultline
\tag 6.22
$$
The first integral is convergent for $\Re s<n/2-k+L$, and is $O(|\xi|^{2\Re s-n+2k-2L})$
and the other terms are meromorphic with simple poles
in the set $n/2+\Bbb N$.  We see that as $|\xi|\to\infty$,
$$
\int_{|\xi|^2}^\infty T^{s-n/2+k}q_k(T,x)\, dT\ \sim\
- \sum_{\ell=0}^{\infty}\frac{q_k^\ell(x)}{s-n/2+k-\ell}\ |\xi|^{2s-n+2k-2\ell},
$$
and this asymptotic formula still holds after analytic continuation.  Now from this we see that
$$
u_{s,k}(x)\ =\ \phi_k(s,x,\xi)\ +\ \psi_k(s,x,\xi)\ =\ \int_0^\infty T^{s-n/2+k}q_k(T,x)\, dT,
$$
which converges for $n/2-k-1<\Re s<n/2-k$,  has an analytic continuation to
$s\in \Bbb C$ with simple poles in $n/2+\Bbb N$, and as $|\xi|\to\infty$,
$$
\phi_k(s,x,\xi)\ \sim\ u_{s,k}(x)\ +\ \sum_{\ell=0}^{\infty}\frac{q_k^\ell(x)}{s-n/2+k-\ell}\ |\xi|^{2s-n+2k-2\ell},
$$
and from this,
$$
\multline
\int_0^1 t^{s-n/2+k}q_k(|\xi|^2t,x)\, dt\ =\ |\xi|^{n-2s-2-2k}\phi_k(s,x,\xi)
\\ \sim\ |\xi|^{n-2s-2-2k}u_{s,k}(x)\ +\ \sum_{\ell=0}^{\infty}\frac{q_k^\ell(x)}{s-n/2+k-\ell}\ |\xi|^{-2\ell-2}.
\endmultline
\tag 6.23
$$
To obtain the asymptotic behavior of $\s(s,x,\xi)$ it remains to compute the asymptotic behavior
of
$$
\int_{0}^1 t^{s-n/2+K}Q_{K}(t,t|\xi|^2,x) \, dt.
$$
We deal with this by making a Taylor expansion of $Q_K(t,T,x)$ as $T\to\infty$:
$$
\multline
\int_0^{1} t^{s-n/2+K} Q_K(t,|\xi|^2t,x)\, dt \\
 =\ \int_0^{1} t^{s-n/2+K}
 \left( Q_K(t,|\xi|^2t,x)-\sum_{\ell=0}^{L-1} |\xi|^{-2\ell-2}Q_K^\ell(t,x) t^{-\ell-1} \right) \, dt
\\ +\ \sum_{\ell=0}^{L-1}\ \left(\int_0^{1} t^{s-n/2+K-\ell-1}
  Q_K^\ell(t,x) \, dt\right)|\xi|^{-2\ell-2}
\endmultline
\tag 6.24
$$
The first integral on the right is $O(|\xi|^{-2L-2})$ provided
$\Re s>n/2-K+L$.  Choosing $K=2L$, we find that it is
$O(|\xi|^{-2L-2})$ on $\Re s>n/2-L$. Putting this together with
(6.23) we find that on the set $\Re s>n/2-L$, as $\xi\to\infty$,
$$
\multline
\s(s,x,\xi)\ \sim\
\sum_{k=0}^{2L-1}|\xi|^{n-2s-2-2k}u_{s,k}(x)\ +\ \sum_{k=0}^{2L-1}
\sum_{\ell=0}^{L-1}\frac{q_k^\ell(x)}{s-n/2+k-\ell}|\xi|^{-2\ell-2}
\\ +\ \sum_{\ell=0}^{L-1}\ \left(\int_0^{1} t^{s-n/2+2L-\ell-1}
  Q_{2L}^\ell(t,x) \, dt\right)|\xi|^{-2\ell-2}
\ +\ O(|\xi|^{-2L-2}).
\endmultline
$$
We note that
$$
Q_{2L}^\ell(t,x)\ =\ \frac1{t^{2L}}\left(q^\ell(t,x)\ -\ \sum_{k=0}^{2L-1}q_k^\ell(x)t^k\right),
$$
and so
$$
\align
\s(s,x,\xi)\ \sim\
\sum_{k=0}^{2L-1} &|\xi|^{n-2s-2-2k}u_{s,k}(x)
\ +\ \sum_{k=0}^{2L-1}\sum_{\ell=0}^{L-1}\frac{q_k^\ell(x)}{s-n/2+k-\ell}|\xi|^{-2\ell-2}
\\ &\qquad +\ \sum_{\ell=0}^{L-1}\ \int_0^{1} t^{s-n/2-\ell-1}
  \left(q^\ell(t,x)\ -\ \sum_{k=0}^{2L-1}q_k^\ell(x)t^k\right) \, dt\ |\xi|^{-2\ell-2}
\ +\ O(|\xi|^{-2L})\\
\ =\
\sum_{k=0}^{2L-1} &|\xi|^{n-2s-2-2k}u_{s,k}(x)
\ +\ \sum_{\ell=0}^{L-1}\ v_{s,k}(x)\ |\xi|^{-2\ell-2}
\ +\ O(|\xi|^{-2L-2})  \tag 6.25\\
\endalign
$$
where $v_{s,k}(x)$ is extended from its domain of convergence by
analytic continuation. It remains to show that all the above
asymptotic formulas actually hold in the sense of symbols.  The
key formulas are (6.22) and (6.24). For (6.24), we set
$$
R_K^L(t,T,x)\ =\ Q_K(t,T,x)-\sum_{\ell=0}^{L-1} Q_K^\ell(t,x) T^{-\ell-1},
$$
so
$$
|\partial_t^k\partial_T^j\partial_x^\al R_K^L(t,T,x)|\ \leq\ CT^{-L-j-1}
$$
and
$$
\partial_x^\al\partial_\xi^\beta\int_0^{1} t^{s-n/2+K}
R_K^L(t,|\xi|^2t,x)\, dt
$$
is a sum of terms of the form
$$
\xi^\gamma\int_0^1  t^{s-n/2+K+j}(\partial_x^\al \partial_T^j R_K^L)(t,|\xi|^2t,x)\, dt
$$
where $|\gamma|=2j-|\beta|$. This is bounded by
$$
|\xi|^{|\gamma|-2L-2j-2}\ \int_0^1 t^{\Re s-n/2+K-L-1} \, dt
$$
which is equal to
$$
C(n,s) |\xi|^{-2L-2-|\beta|}
$$
provided $\Re s>n/2-K+L$.  This shows that
$$
\int_0^{1} t^{s-n/2+2L} R_{2L}^L(t,|\xi|^2t,x)\, dt
$$
is a symbol of order  $-2L-2$ when $\Re s>n/2-L$.

Similarly, for (6.22), set
$$
r_k^L(T,x)\ =\ q_k(T,x)-\sum_{\ell=0}^{L-1} q_k^\ell(x)T^{-L-1},
$$
so
$$
|\partial_x^\al\partial_T^j r_k^L(T,x)|\ \leq\ CT^{-L-j-1}.
$$
Then
$$
\partial_{\xi_m} \int_{|\xi|^2}^\infty T^{s-n/2+k} r_k^L(T,x)\, dT
\ =\ -2\xi_m |\xi|^{2s-n+2k}r_k^L(|\xi|^2,x)
$$
is a symbol of order $2\Re s-n+2k-2L-1$, and
$$
\left|\partial_{x}^\al \int_{|\xi|^2}^\infty T^{s-n/2+k} r_k^L(T,x)\, dT\right|\ \leq\
C\,|\xi|^{2\Re s-n+2k-2L}.
$$
It follows that
$$
\phi_k(s,x,\xi)\ -\ u_{s,k}(x)\ +\ \sum_{\ell=0}^{L-1}q_k^\ell(x)|\xi|^{2s-n+2k-2\ell}
\ =\ \int_{|\xi|^2}^\infty T^{s-n/2+k} r_k^L(T,x)\, dT
$$
is a symbol of order $2\Re s-n+2k-2L$, and  putting this together,
the $O(|\xi|^{-2L-2})$ error in (6.25) is actually a symbol of
order $-2L-2$. This completes the proof of Lemma 6.5.\qed
\enddemo

\medskip
\medskip
\medskip



\vskip20pt

\centerline{References}

\roster

\item"{[Br]}"
T. Branson,
Sharp inequalities, the functional determinant, and the
complementary series,
{\it Preprint}.

\item"{[BCY]}"
T. Branson, S. Y. A. Chang and P. Yang,
Estimates  and extremals for zeta function determinants on four-manifolds,
{\it Commun. Math. Phys.} {\bf 149}, (1992), 241-262.

\item"{[CGY]}"  S. Y. A. Chang, M. Gursky and P. Yang,
An equation of Monge-Ampere type in
conformal geometry, and four-manifolds of positive Ricci curvature,''
{\it Annals of Math.}, {\bf 155} (2002) pp. 1-79.

\item"{[CQ]}"
S. Y. A. Chang  and J. Qing,
Zeta Function Determinants on Manifolds with Boundary,
{\it Math. Research Letters} {\bf 3}, 1-17  (1996).

\item"{[CY]}"
S. Y. A. Chang and P. Yang,
Extremal Metrics of zeta function determinants on 4-manifolds,
{\it Annals of Math.} {\bf 142} (1995) 171-212.

\item"[Ch]"
I. Chavel,
{\it Eigenvalues in Riemannian geometry},
Pure and Applied Math. {\bf 115},
Academic Press, (1984).

\item"[Gi]"
P. Gilkey,
{\it The index theorem and the heat equation.}
Publish or
Perish, Inc., Boston, Mass., (1974).

\item"[GS]"  V. Guillemin and S. Sternberg,   Some remarks on I.
M. Gelfand's work,  {\it Izrail M. Gelfand,  Collected Papers},
{\bf Vol. I}, Springer-Verlag, Berlin, (1987), 831--836.

\item"[H\"{o}]"
L. H\"{o}rmander,
The spectral function of an elliptic operator.
{\it Acta Math.} {\bf 121}, 193-218 (1968).

\item"[Mo1]"
C. Morpurgo,
Sharp trace inequalities for intertwining operators on $S\sp n$ and $R\sp n$.
{\it Internat. Math. Res. Notices}  (1999),
{\bf  20}, 1101--1117.

\item"[Mo2]"
C. Morpurgo,
Sharp inequalities for functional integrals and
traces of conformally invariant operators,
{\it Duke Math. Journal}, {\it to appear}.

\item"[Ok1]"
K. Okikiolu,
Critical metrics for the determinant of the Laplacian in odd dimensions.
{\it Ann. of Math.}
 {\bf 153} (2001),  471--531.

\item"[Ok2]"
K. Okikiolu, Critical metrics for the zeta function of  a family
of scalar Laplacians, {\it Preprint.}

\item"[OW]"
K. Okikiolu and C. Wang,  Hessian of the zeta function for the
Laplacian on forms. {\it Preprint.}

\item"{[OPS1]}"
B. Osgood, R. Phillips and P. Sarnak,
Extremals of Determinants of Laplacians,
{\it  J.  Funct. Anal.} {\bf  80} (1988), 148-211.

\item"{[OPS2]}"
B. Osgood, R. Phillips and P. Sarnak,
Compact isospectral sets of surfaces,
{\it  J.  Funct. Anal.}  {\bf 80} (1988), 212-234.

\item"{[OPS3]}"
B. Osgood, R. Phillips and P. Sarnak,
Moduli space, heights and isospectral sets of plane domains,
{\it  Annals of  Math.}  {\bf 129} (1989), 293-362.

\item
"[{ Ri}]"
{K. Richardson},
{Critical points of the determinant of the Laplace operator},
{\it Jour.  Funct. Anal.}, {\bf 122}, 52-83 (1994).

\item
"[Se]" R. Seeley,
Complex Powers of an Elliptic Operator,
{\it Proc.  Symp. on Singular Integrals},
AMS {\bf 10}, (1967), 131-160.

\item
"[Tr]"
F. Treves,
{\it Introduction to pseudodifferential and Fourier integral operators.}
{\bf Vol. 1.},
Plenum Press,  (1980).

\endroster
\medskip
\medskip
\medskip

\centerline{Kate Okikiolu}
\centerline{University of California, San Diego}
\centerline{okikiolu\@math.ucsd.edu}

\end